\documentclass[10pt]{article}

\usepackage{fancyheadings}
\renewcommand{\headrulewidth}{0.5pt}
\usepackage[all]{xy}
\usepackage[]{amsmath,amscd,amsfonts,amsthm,amssymb}
\usepackage{amssymb}
\input{diagxy}
\usepackage{times}
\usepackage{epic}
\usepackage{graphicx}
\usepackage{graphics}
\usepackage{cjhebrew}
\usepackage{psfrag}
\usepackage[lflt]{floatflt}
\title{\LARGE \sc Zeno's Paradoxes. A Cardinal Problem\\ 
\vspace{0.5cm}\large I. On
Zenonian Plurality}
\author{{\em Karin Verelst}}
\date{}

\newfam\msbfam
\font\tenmsb=msbm10                     \textfont\msbfam=\tenmsb
\font\sevenmsb=msbm7            \scriptfont\msbfam=\sevenmsb
\font\fivemsb=msbm5
\scriptscriptfont\msbfam=\fivemsb

\newcommand{\pw}{\footnotesize}
\newcommand{\bq}{\begin{quote}}
\newcommand{\eq}{\end{quote}}

\newcommand{\pijll}{\longleftarrow}
\newcommand{\pijlr}{\longrightarrow}
\newcommand{\pijltje}{\rightarrow}

\newcommand{\aequi}{\Leftrightarrow}
\newcommand{\lppijl}{\longleftrightarrow}

\newcommand{\krulpijlr}{\hookrightarrow}
\newcommand{\order}{\leqslant}

\newcommand{\M}{\mid}
\newtheorem{Th}{Theorem}
\newtheorem{ax}{Axiom}
\newtheorem{lm}{Lemma}
\newtheorem{df}{Definition}
\newtheorem{pr}{Proposition}
\newtheorem{cl}{Corollary}
\newtheorem{re}{Remark}
\newtheorem{as}{Assumption}
\newtheorem{wg}{Wild Guess}
\newtheorem{ex}{Example}
\newcommand{\bth}{\begin{Th}\hspace{-5pt}{\bf .} \ }
\newcommand{\Eth}{\end{Th}}
\newcommand{\bax}{\begin{ax}\hspace{-5pt}{\bf .} \ }
\newcommand{\eax}{\end{ax}}
\newcommand{\blm}{\begin{lm}\hspace{-5pt}{\bf .} \ }
\newcommand{\elm}{\end{lm}}
\newcommand{\bdf}{\begin{df}\hspace{-5pt}{\bf .} \ }
\newcommand{\edf}{\end{df}}
\newcommand{\bpr}{\begin{pr}\hspace{-5pt}{\bf .} \ }
\newcommand{\epr}{\end{pr}}
\newcommand{\bcl}{\begin{cl}\hspace{-5pt}{\bf .} \ }
\newcommand{\ecl}{\end{cl}}
\newcommand{\bre}{\begin{re}\hspace{-5pt}{\bf .} \ }
\newcommand{\ere}{\end{re}}
\newcommand{\bas}{\begin{as}\hspace{-5pt}{\bf .} \ }
\newcommand{\eas}{\end{as}}
\newcommand{\bwg}{\begin{wg}\hspace{-5pt}{\bf .} \ }
\newcommand{\ewg}{\end{wg}}
\newcommand{\bex}{\begin{ex}\hspace{-5pt}{\bf .} \ }
\newcommand{\eex}{\end{ex}}

\newcommand{\bit}{\begin{itemize}}
\newcommand{\eit}{\end{itemize}\par\noindent}
\newcommand{\ben}{\begin{enumerate}}
\newcommand{\een}{\end{enumerate}\par\noindent}
\newcommand{\beq}{\begin{equation}}
\newcommand{\eeq}{\end{equation}}
\newcommand{\beqa}{\begin{eqnarray*}}
\newcommand{\eeqa}{\end{eqnarray*}\par\noindent}
\newcommand{\beqn}{\begin{eqnarray}}
\newcommand{\eeqn}{\end{eqnarray}\par\noindent}
\def\nn{\mbox{I\hspace{-0.5mm}N}} 

\def\rr{\mbox{I\hspace{-0.5mm}R}}

\def\bb{\mbox{I\hspace{-0.5mm}B}}

\def\pp{\mbox{I\hspace{-0.5mm}P}}

\def\wit{\hspace{1pt}}
\def\isub#1#2{\mbox{$\underset{\hspace{#2pt}\iota}{#1}$}}
\def\csup#1#2{\mbox{$#1\hspace{-9pt}\hspace{-#2pt}^{\mbox{
				\tiny $^{\smallfrown}$}}$}}
\def\spil#1#2{\mbox{$#1\hspace{-4pt}\hspace{-#2pt}^{\mbox{\footnotesize 
'}}$\wit}}
\def\spir#1#2{\mbox{$#1\hspace{-4pt}\hspace{-#2pt}^{\mbox{\footnotesize 
`}}$\wit}}
\def\acca#1#2{\mbox{$#1\hspace{-4pt}\hspace{-#2pt}^{\mbox{\tiny
$^{\prime}$}}$\wit}}
\def\accg#1#2{\mbox{$#1\hspace{-4pt}\hspace{-#2pt}^{\mbox{\tiny
$^{\backprime}$}}$\wit}}
\def\spilacca#1#2{\mbox{$#1\hspace{-8.5pt}\hspace{-#2pt}^{\mbox{
			\footnotesize '\tiny $^{\prime}$}}$}}

\def\spiracca#1#2{\mbox{$#1\hspace{-8.5pt}\hspace{-#2pt}^{\mbox{
			\footnotesize `\tiny $^{\prime}$}}$}}

\def\csupspil#1#2{\mbox{$#1\hspace{-4pt}\hspace{-#2pt}^{\csup{
			\mbox{\footnotesize '}}{-2}}$}}

\usepackage[]{amsmath,amscd,amsfonts,amsthm,amssymb}
\usepackage{fancybox}

\newenvironment{myboxqd}
           {\begin{Sbox}\begin{minipage}{1\textwidth}}
           {\end{minipage}\end{Sbox}\Ovalbox{\TheSbox}}

\begin{document}

\maketitle
\vspace{-0.7cm}
{\small
\centerline{FUND-CLEA}
}
{\scriptsize
\centerline{Vrije Universiteit Brussel}
\centerline{kverelst@vub.ac.be}
}

\thispagestyle{plain}
\cfoot{}

\bigskip
\bigskip

\begin{flushright}{\small {\em Du pr\'{e}sent, rien d'autre que du 
pr\'{e}sent\\ {\em
Sartre, La naus\'{e}e}}} \\
\end{flushright}

\bigskip

\frenchspacing {\sc Introduction}.--- It will be shown in this 
contribution that the
Received View on Zeno's paradoxical arguments is untenable. Upon a 
close analysis of
the Greek sources\footnote{The reference textcritical edition for the 
fragments [B] and
related testimonia [A] is: H. Diels and W. Kranz [DK in what follows. 
See the list of sigla at the end
of this paper], {\em Fragmente der Vorsokratiker}, to the numbering 
of which I will
comply in accordance with scholarly tradition.}, it is possible to do 
justice to
Simplicius's widely neglected testimony, where he states: {\em In his 
book, in which many
arguments are put forward, he shows {\bf
\em in each} that a man who says that there is a {\bf \em plurality} 
is stating something
{\bf \em contradictory}} [{\small DK} {\scriptsize{29B 2}}]. Thus we 
will demonstrate
that an underlying structure common to both the Paradoxes of 
Plurality (PP) and the
Paradoxes of Motion (PM) shores up {\em all}
\/his arguments.\footnote{In this we reckon in Owen a precursor, although our
analysis of Zeno's arguments will be very different from his. See 
G.E.L. Owen, ``Zeno
and the  mathematicians'', {\em Proceedings of the Aristotelian
Society}, {\bf 8}, 1957.} This structure
bears on a correct --- Zenonian --- interpretation of the concept of 
``division through
and through'', which takes into account \vspace{-0.5 mm} the often 
misunderstood
Parmenidean legacy, summed up concisely in the deictical dictum
$\tau\accg{o}{0} \ \spil{\epsilon}{-1}\accg{o}{-1}\nu
\ \spilacca{\epsilon}{-2}\sigma\tau\iota$ : ``the Being-Now 
is''.\footnote{K. Riezler,
{\em Parmenides. Text,
\'{U}bersetzung, Einf\"{u}hrung und Interpretation}, Vittorio 
Klostermann, Frankfurt,
1970, p. 45-50. Parmenides's $\tau\accg{o}{0} \ 
\spil{\epsilon}{-1}\accg{o}{-1}\nu
\ \spilacca{\epsilon}{-2}\sigma\tau\iota$ {\em [to eon esti]} should 
indeed be translated
as {\em The Now-Being is}. In the present everything {\em
is}. That the in origin dialectal difference between $\tau\acca{o}{-2} \
\spilacca{o}{-1}\nu$ and $\tau\accg{o}{0} \ 
\spil{\epsilon}{-1}\accg{o}{-1}\nu$ had
acquired philosophical significance becomes explicit in Diogenes of 
Appolonia, a
contemporary to Zeno, where
$\tau\accg{\alpha}{-2} \ \spilacca{o}{-1}\nu\tau\alpha$ {\em [ta 
onta; the beings]}
are stable essences, while
$\tau\accg{\alpha}{-1} \
\spil{\epsilon}{-1}\acca{o}{-1}\nu\tau\alpha \ 
\nu\csup{\upsilon}{-2}\nu$ {\em [ta
eonta nun; the beings-now]} are instable phaenomenological things. See L.
Couloubaritsis, {\em La Physique d'Aristote}, Ousia, Bruxelles, 1997, 
p. 308.} The
feature, generally overlooked but a key to a correct understanding of 
all the arguments
based on Zeno's divisional procedure, is that {\em they do not 
presuppose space, nor
time.} Division merely requires extension of an object present to the 
senses and takes
place simultaneously. This holds true for both PP and PM! Another 
feature in need of
rehabilitation is Zeno's plainly avered but by others blatantly 
denied phaenomenological
--- better: {\em deictical} --- realism. Zeno  nowhere denies the 
reality of plurality and
change. Zeno's arguments are not a {\em reductio}, if only because 
the logical prejudice
that something which implies paradoxes cannot ``really be there'' is 
itself still
unthinkable, since hypothetical thinking does not yet exist. When one 
speaks, one does so
not about a possible world, but about {\em this} world:
$\kappa\acca{o}{-1}\sigma\mu o\varsigma \ \tau\acca{o}{-1}\delta\epsilon$ \/ as
Heraclitus calls her [{\small DK} {\scriptsize{22B 
30}}].\footnote{Conche translates ``ce
monde-ci''. M. Conche, {\em H\'{e}raclite. Fragments}, Presses 
Universitaires de France,
1986/1998, pp. 279-280.} Zeno merely shows that, when someone states 
plurality, he
inevitably states a contradiction, exactly as Simplicius claims. In what are
traditionally considered the plurality arguments, this contradiction 
appears in the
stature of
$\mu\acca{\epsilon}{-1}\gamma\alpha\lambda\alpha \ 
\kappa\alpha\accg{\iota}{-1} \
\mu\iota\kappa\rho\accg{\alpha}{-1}$, the large[s]-and-small[s] [{\small DK}
{\scriptsize{29B 1}}]. After subsuming PM under the simultaneous 
divisional procedure
proper to PP, it will be indicated how the received view on the 
former can easely be
derived by the introduction of time as a (non-Zenonian) praemiss, 
thus causing their
collaps into arguments which can be approached and refuted by 
Aristotle's limit-like
concept of the ``potentially infinite'', which remained --- in 
different disguises --- at
the core of the refutational strategies that have been in use up to 
the present. A
mathematical representation will be given for Zeno's simultaneous 
divisional procedure
which fully reckons Aristotle's dictum, revealingly enough never 
discussed in relation
to Zeno, where he says: {\em For in two ways it can be said that a 
distance or a period or
any other continuum is infinite, viz., with respect to the partitions 
or with respect to
the parts} [{\em Phys.} Z, 2, 263a (24-26)].\\

{\sc Zeno's Infinite Division and the Paradoxes of Plurality}.--- I 
said that the sources
testify for the unity of Zeno's arguments, in that they {\sc all} are 
arguments on
plurality. In this paragraph we are going to read in some detail the 
fragments which are
seen traditionally as the paradoxes of plurality, in order to 
understand how Zeno gets to
his plurality and what he has in mind when he uses that term. A comment by
Simplicius will serve as our guideline:
$\kappa\alpha\tau\accg{\alpha}{-0.5} \
\tau\accg{o}{-1} \ \pi\lambda\csup{\eta}{-1}\theta o\varsigma \
\spilacca{\alpha}{-1}\pi\epsilon\iota\rho o\nu \ \spil{\epsilon}{-1}\kappa \
\tau\csup{\eta}{-1}\varsigma \ \delta\iota \chi o\tau
o\mu\acca{\iota}{-1.5}\alpha\varsigma \
\spilacca{\epsilon}{-1}\delta\epsilon\iota\xi\epsilon \,$ {\em [Thus] 
he demonstrated
numerical infinity by means of dichotomy.} [{\em Phys.}, 140 (27)].

\noindent Three observations are in order here. First, the verbal 
form of the main clause
is not neutral with regard to Zeno's achievement: the conjugation used is
$\spilacca{\epsilon}{-1}\delta\epsilon\iota\xi\epsilon \,$, the 
($3^{rd}$ sing.) aorist
of $\delta\epsilon\iota\kappa\nu\acca{\upsilon}{-1}\nu\alpha\iota$ 
\,{\em [deiknunai]}
(`to point at, to indicate'; here in the sense of `to show',`to 
demonstrate'), so as to
mark out unambiguously that Zeno demonstrated this once and for all, 
exactly as when we
speak of ``G\"{o}del's proof'' as an acquired, definite result. 
Second, the word {\em
apeiron}, customary translated as `infinity', occurs in Zeno's own 
tekst and means
literally `unbounded'. It derives from the archaic {\em a-peirar};
with $\pi\epsilon\csup{\iota}{-1}\hspace{-0.5 mm}\rho\alpha\rho$ : 
rope, knot or
bond\footnote{R.B. Onians, {\em The Origins of European
Thought}, Cambridge University Press, Cambridge, 1951[1994], p. 
314-317; 332 sq.}, i.e.,
something that has to be put around or upon something else from the 
outside, not just an
end or a limit which is intrinsic to it. It is rather like a fence 
enclosing a meadow.
{\em Okeanos}, the primal sea, is said to be {\em apeiron}; the image 
here being that of
a lonely ship (the Earth) in the middle of a sea with no land in 
sight. So when used to
qualify an abstract entity or process, the connotation will be `not 
stopped by somebody
or something, uninterrupted'. One will be reminded of the `bonds of 
Necessity' invoked by
Parmenides in his Poem [{\small DK} {\scriptsize{28B 8 (30)}}]. Third 
and last, the word
$\delta\iota\chi o\tau o\mu\acca{\iota}{-1.5}\hspace{-0.5 mm}\alpha$ 
{\em [dichotomy]}
does {\em not} occur in Zeno's own text, and we know that this is not 
for want of
relevant text pieces. How then do we know that Zeno has in view such 
a division? We will
see in a minute why the {\em procedure} Zeno describes in his text 
cannot be anything
else than a `division through and through'. But there are other, 
related testimonies. A
fragment ascribed by Porphyry to Parmenides does mention
$\delta\iota\alpha\acca{\iota}{-1}\hspace{-1.5 mm}\rho\epsilon\tau 
o\nu$ {\em [division]}
explicitly. Simplicius [{\em Phys.}, 140 (21)] points out that this 
attribution must be
fallacious: {\em For no such arguments figure among the Parmenidean 
[texts] and the
majority of our information refers the difficulty from dichotomy to 
Zeno}. This is
confirmed by Philoponus [{\em In Physica} 80(26-27)].\footnote{LEE, 
p. 22. S. Makin,
``Zeno on Plurality'' {\em Phronesis}, {\bf 27}, 1982, pp. 223-238, 
gives Themistius's
commentary on the {\em Physics} as further evidence. I follow Vlastos 
for the emendation
of the word ``texts'' in the fragment quoted; VLAS, p. 231. But even 
if the fragment were
Parmenidean, it would not lose its relevance, given the close 
doctrinal relationship
between the two men; so W.E. Abraham, ``The nature of Zeno's Argument 
Against Plurality in
[{\small DK} {\scriptsize{29B 1}}]'', {\em Phronesis}, {\bf, 17}, 
1972, pp. 40-52.} In
the Porphyry text it is stipulated: {\em since it is alike throughout
$\pi\acca{\alpha}{-1}\nu\tau\isub{\eta}{0}$ [pant\={e}i],\vspace{-2mm} if it is
divisible, it will be divisible throughout alike, not just here but not
there}.\footnote{VLAS, p. 229; LEE, pp. 12, 20-23. A close parallel 
to this argument can
be found in Aristotle's book on becoming, [{\em De gen. et cor.}, I. 
2, 316a16 sq. and
325 a8].} The relevance of the words `here'`and `there' cannot be 
overestimated, for in
Zenonian terms they serve to deictically define extension. Simplicius 
moreover adds: {\em
There is no need to labour the point; for such an argument is to be 
found in Zeno's own
book. (...) Zeno writes the following words (...)}\footnote{See LEE, 
p. 21 (the source of
the translation). Compare DK p. 257, ftn. 5 with [{\small DK} 
{\scriptsize{28B 8
(22)}}].}; after which follows a literal quotation:\\

The second Paradox of Plurality (infinite divisibility) [{\small DK} 
{\scriptsize{29B
3}}]\footnote{Concerning the translation of Zeno's arguments: I made 
use of LEE, KRS and
other sources to be mentioned in case, but nowhere I follow them 
completely, sometimes
--- I admit --- to the detriment of the English used. This is because 
I chose to contract
O'Flaherty's methodological advise (in her book on sexual metaphor in 
Ancient Indian
mythology) as completely as possible: {\em In the  first analysis, it 
pays to be
literal-minded}. W.D. O'Flaherty, {\em Women, Androgynes and other 
mythical Beasts},
University of Chicago Press, Chicago, 1980, p. 5. I can only hope the 
reader will be
indulgent with respect to this choice.}

\begin{quotation} [Simplicius, {\em In Aristotelis physicorum}, 140 
(27)] \, {\em For in
his proof that, if there is plurality, the same things are both 
finite and infinite, Zeno
writes the following words: ``if they are many [things], they by 
necessity are as many as
they are, not more nor less. But if they are as many as they are, 
they will be finite
[bounded, peperasmena]. But if they are many, they will be infinite 
[unbounded, apeiron].
For there will always [aei] be others [hetera] in between [metaxu] of 
the beings, and
there again others in between.'' Thus he demonstrated infinity by 
means of dichotomy.}
\end{quotation}

\noindent It is remarkable that Simplicius is so firm in his 
statement, while we would be
tempted to find the argument at first glance rather weak. Remember Simplicius's
introduction: {\em he shows in each [argument] that a man who says 
there is a plurality
is stating something contradictory}. The contradiction apparently is 
that if ``they are
many, they will be both bounded and unbounded''.  Whence does this 
contradiction arise?
The point is clearly somewhere in the sentence: {\em For there will 
always be others in
between of the beings, and there again others in between}. Even if 
not mentioned
explicitly, this cannot be other than some kind of divisional 
procedure, exactly as we
conceive of fractions to mentally break a line. Zeno's 
plurality-argument constitutes
{\bf the first Gedanken-experiment} in scientific history! It 
furthermore states in an
accurate way that the number of parts obtained should at least be 
{\em dense}, in the
mathematical sense of that word.\footnote{I am not ``reading this 
into'' Zeno; his
formulation is by far the closest you can get to it in words. Let's 
make the point by
comparing him to a standard textbook definition: {\em Let $a$ and $b$ 
be two real numbers
with $a < b$. We can always find a real number $x$ between $a$ and 
$b$.} See K.G.
Binmore, {\em The Foundations of Analysis: A straightforward 
Introduction}, Cambridge
University Press, Cambridge, 1980, p. 74. The point is so obvious 
that one wonders why it
is not made more often.} To be sure, this also means that atoms as 
straightforward
``least parts'' are non-Zenonian. It makes moreover plain that `being 
(un)bounded' really
means `being (un)limited in number'. I am running a bit ahead of my 
argument when I say
that this implies as well that they are at least countably infinite. 
We already know that
  limitation comes about by something external, by an obstacle or by 
interruption of an
ongoing process. It inevitably follows that Zeno had in mind a 
division which is both
symmetrical and nondirected, i.e. one in which {\em all} parts 
undergo the {\em same}
divisional process.\footnote{However, there IS an unexpressed 
hypothesis crucial to Zeno's
procedure, namely that division be one-dimensional, and in 
imagination presented as
horizontal. This allows him to go over smoothly form PP to PM, but it 
has non-innocent
mathematical consequences, on which we will come back. The relevance 
of this for the
mathematics involved was pointed out to me by Bob Coecke (Oxford).} 
Zeno's terminology
speaks for itself, as both
$\mu\epsilon\tau\alpha\xi\acca{\upsilon}{0}$\, {\em [metaxu]} 
`amidst', `in between equal
parts or things' and $\spiracca{\epsilon}{-2}\tau\epsilon\rho\alpha$\ 
{\em [hetera]}
`others', where the singular {\em heteron} indicates `the other of 
two', testify. With
regard to the latter Diels-Kranz put it plain and simple in their apparatus:
$\spiracca{\epsilon}{-2}\tau\epsilon\rho o\nu$  {\em : 
Dichotomie!}\/\footnote{DK, vol.
I, p.
255.} 

Thus for anyone sharing the terminological sensitivity that goes along with having
ancient Greek as a mothertongue, Zeno's intentions were immediately 
clear. But we do
still not know why applying this procedure would amount into 
paradoxical results. To
answer this question, we have to extract more information from 
another variant of the
argument, known as the `first' paradox of plurality, because in it 
the intended situation
is even more explicitly exposed:\\

The first Paradox of Plurality (finite extension) [{\small DK} 
{\scriptsize{29B 1
\& B 2}}]

\begin{quotation} [Simpl., {\em Phys}, 140 (34)] \, {\em The infinity 
of magnitude he
showed previously by the same reasoning; for, having first shown that 
``if a being had no
magnitude, it would not be at all'', he proceeds ``but if it is, then 
each one must
necessarily have some magnitude [megethos] and thickness and keep one 
away [apechein
heteron] from the other [apo heteron]. And the same reasoning holds 
for any [part]
jutting out [peri to prouchontos]; for this too will have extended 
magnitude and jut out.
But to say this once is as good as saying it forever [aei]; for none 
will be the last
[eschaton] nor the one [heteron] will be unrelated to another one 
[pros heteron]. So if
it is a plurality, it by necessity will be \vspace{-1.5 mm} many 
small [ones] and many
large [ones] {\bf \em [mikra kai megala]} --- 
$\mu\iota\kappa\rho\accg{\alpha}{-1} \
\tau\epsilon \ \epsilon\csupspil{\iota}{0}\hspace{-1.5 mm}\nu\alpha\iota
\ \kappa\alpha\accg{\iota}{-1} \ 
\mu\acca{\epsilon}{-1}\gamma\alpha\lambda\alpha$ --- ;
so many small[s] as to have no magnitude [m\={e} megethos], so many 
large[s] as to be
unbounded [apeiron]''.}
\end{quotation}

\begin{quotation} [Simpl., {\em Phys}, 139 (5)] \, {\em In one of 
these arguments he
shows that if there is plurality [polla], then it is both many large 
[ones] and many
small [ones]\footnote{I know of no translations which renders this 
sentence coorectly. The
verb used is \/$\spil{\epsilon}{-1}\sigma\tau\acca{\iota}{-2}$, ``it 
is'', while {\em
megala} and {\em mikra} obviously are plurals. This moreover rules 
out the translation
``they must be both small and large'' for the last sentence in [{\small DK}
{\scriptsize{29B 1}}], as for instance in LEE.}, so many large[s] 
[megala] as be infinite
[apeiron] in magnitude, so many small[s] [mikra] as to have no 
magnitude at all. In this
same argument he shows that what is without magnitude, thickness and 
bulk would not be at
all. ``For'', he says, ``if it were added to some other being, it 
would not make it
bigger; because being of no magnitude, when added, it cannot possibly 
make it grow in
magnitude. And thus the added would in fact be nothing. So if, when 
taken away, the other
[being] will not be any less, and again will not, when added, 
increase, then it is clear
that the added and then again taken away was nothing.''}\footnote{The 
hesitation in the
standard translations with respect to the last line seems unnecessary 
when one takes one
of the very few Aristotelian texts into account that deal explicitly 
with plurality ---
interestingly enough in [{\em Met.} 1001b7-19], i.e., not in the {\em 
Physics} --- and
where an almost literal quotation of Zeno's words is present: {\em 
For, he says, that
which makes [something else] no larger, when added, and no smaller, 
when subtracted, is
not an existent}. Translation with [ ]: VLAS, p. 238. }
\end{quotation}

\noindent Our aim is to try to understand Zeno in the way he 
understood himself, so we
will follow Zeno's own argumentation as closely as possible and try 
to avoid any
assumption imposed on him by later times. This explains my a bit 
awkward translation: the
precise use of singular verb forms and plural substantives in the 
description is essential
to a correct understanding, and should as much as possible be 
preserved. We already found
that Zeno must have had in mind some kind of division. From the first 
lines of [{\small
DK} {\scriptsize{29B 1}}] it is clear that we are dealing with an 
object that has material
extension; the verb `proechein' is {\em normally said of contiguous 
parts, one of which
is thought of as sticking out from or extending beyond the 
other}.\footnote{VLAS,
``Plurality'', p. 226.} Plurality is coined here in terms of the 
relation between parts
and whole in an extended object with a `here' and a `there': {\em The 
discrimination of
any two such parts in any existent I shall refer to as a 
``division''}.\footnote{VLAS,
``Plurality'', p. 226.} The deictic ``first person''\footnote{E. 
Benv\'{e}niste, ``Le langage et
l'exp\'{e}rience humaine'', in: {\em Probl\`{e}mes de linguistique 
g\'{e}n\'{e}rale II},
Gallimard,  Paris, 1966, p. 69.} standpoint manifest in the lacking 
of temporality in
Zeno's {\em present}-ation apparently has a spatial counterpart. 
``Space'' nor ``time''
exist as independent backgrounds against wich a mentally 
representable event takes
place\footnote{J. Bollack, H. Wismann, {\em H\'{e}raclite ou la 
s\'{e}paration}, Editions
de Minuit, Paris, 1972,  p. 49.}; it is by indication that
the validity of the argument is {\em shown}. Let us now see, on the 
basis of these two
fragments, what characterises this division, and join our conclusions 
with those attained
on the first fragment. There has been considerable controversy in the 
literature on this
subject. It is appropriate to follow here Abraham's terminological 
distinctions between
bipartite and tripartite division, and between simple `division at 
infinity' and
`division throughout', i.e., stepwise applied to the last product of 
division or applied
to all obtained parts equally. The bipartite/tripartite controversy was settled
(methinks convincingly) by Vlastos. He argues for the symmetrical 
variant, offering,
apart from considerations on the impact of symmetrical proportion on 
the archaic mind, a
number of linguistic arguments, which strengthen the few we offered 
with regard to
[{\small DK} {\scriptsize{29B 1}}]. The image he presents of Zeno's 
procedure is the one
now generally accepted: imagine a rod, divide it in two equal parts, 
take the right hand
part, divide it the same way, and repeate this procedure {\em ad 
infinitum}. One then
obtains the physical aequivalent of a mathematical sequence of 
geometrically decreasing
parts. This is exemplified in Vlastos's translation of the clause
$\pi\epsilon\rho\accg{\iota}{-1} \ \tau o\csup{\upsilon}{-1} \ \pi\rho
o\acca{\upsilon}{-1}\chi o\nu\tau o\varsigma$ is ``[And the same 
reasoning applies] to
{\em the} projecting [part]'', the part which remains to be further 
divided, that is.

Zeno's conclusion inevitably is that ``a finite thing is infinite'', not a
paradoxcical, but a ridicoulous statement. The assertion then mostly 
follows that Zeno's
reasoning is based on a lack of mathematical knowledge for it is 
evident, isn't it, that
the sum of an infinite series can very well have a finite total. 
Vlastos, who is familiar
with the subtlety and profundity of Ancient Greek thought  has the 
elegance at least to
look for other explanations of Zeno's supposed blatant 
error.\footnote{VLAS, p. 234.} It
is nevertheless true that on this interpretation the unity of Zeno's 
arguments is
respected at least to the extend that two of the four Zenonian 
paradoxes of motion can be
subsumed under the same model, in the received view on PM. Its major 
flaw with respect to
at least PP is, however, that it cannot be upheld, because, 1) the 
density-property
explicitly mentioned in [{\small DK} {\scriptsize{29B 3}}] remains 
unexplained; 2) the
Porphyry text gives terminological evidence in favour of a `division 
throughout': {\em
since it is everywhere [pant\={e}i] homogeneous, if it is divisible, 
it will be divisible
everywhere [pant\={e}i] alike}\footnote{The
\vspace{-1 mm} Greek is
$\pi\acca{\alpha}{-1}\nu\tau\isub{\eta}{0}$ {\em [pant\={e}i]}: 
overall, everywhere,
carries the idea of an undiscriminated application.}; and 3) an 
essential part of the
claim stated in [{\small DK} {\scriptsize{29B 1}}] is not taken into 
account. Therein it
is said that {\em if it is a plurality, it by necessity will be many 
small [ones] and
large [ones]; so many small[s] as to have no magnitude, so many 
large[s] as to be
infinite}. On the interpretation discussed up to now the first part 
of the assertion is
plainly neglected. The idea of many commentators is in all likelihood 
that from the quote
cited by Simplicius at the beginning of the same fragment [{\em if 
what is had no
magnitude, it would not exist at all}], it can be inferred that one 
can simply dismiss
this possibility. A more subtle view of the matter credits Zeno with 
the contemplation of
things being but without magnitude, like unextended mathematical 
points.\footnote{The
French author P. Tannery introduced in the modern literature the idea that set
theory and the paradoxes appearing in it should be related to the 
work of Zeno. P.
Tannery, ``Le concept scientifique du continu. Z\'{e}non d'El\'{e}e 
et Georg Cantor",
{\em Revue philosophique de la France et de l'\'{e}tranger}, {\bf 
20}, 1885, p. 397 sq.}
The reasoning supposedly goes as follows: Zeno says that a) infinite 
division leads to an
infinite number of final, indivisible parts which still do have a 
magitude, because, b)
if they would not have magnitude, they would not be at all, and so 
the object of which
they are part would not be at all. And an infinite number of parts 
possessing, however
small, finite magnitude, would give us an object infinitely big. But 
c) given that the
division is complete (`throughout'), no parts with finite magnitude can
remain\footnote{An assumption which is made explicit is the 
``Porphyry text'': {\em if
any part of it is left over, it has not yet been divided throughout 
[pant\={e}i]}. VLAS,
p. 229.}, therefore d) the object constituted by them will be 
infinitely small. Thus,
upon Zeno's argument, a finite thing consisting of a plurality would 
be either infinitely
large or infinitely small. The subtle variant of the standard 
interpretation thus offers
us Zeno's argument as a dilemma. Owen develops the dilemma explicitly 
in terms of
divisibility.\footnote{G.E.L. Owen, {\em
op. cit.}.} This is also what Huggett does, by deriving as horns from the
proposition {\em The points have either zero length or finite 
length}\, the conclusions
{\em C1. The total length of the segment is zero} versus {\em C2. The 
total length of the
segment is infinite}.\footnote{N. Huggett, {\em o.c.}, pp. 44-45.} 
This viewpoint
necessitates anyhow another interpretation of his infinite division, 
namely that {\em
every} resulting part will be subject to the infinite series of 
stepwise divisions. One
then obtains a countably infinite number of decreasing sequences, 
resulting in a
countably infinite number of dimensionless endpoints. This reading of 
the argument has a
venerable tradition and it respects the historical order of things, 
since this is the way
the atomists interpreted it.\footnote{Epicurus in his Letter to 
Herodotus, 56-57. M.
Conche, {\em Epicure, Lettres et Maximes}, PUF, Paris, 1987/1999, pp. 
108-111. See also
the commentary on pp. 147-151.} Alas, Zeno nowhere presupposes 
material atoms. This
process is taken to be carried through somehow up to the moment the 
{\em unextended}
points composing the (rational part of) real line are reached, the 
``infinitieth''
element of every sequence. It is eventually concluded that Zeno 
commits a fatal, yet this
time logical, fallacy. Indeed, on the assumption --- necessary if the 
enquiry into the
consequences of the concept `plurality' is to be exhaustive --- that 
such an ordinal
[i.e., stepwise] infinite process could be {\em completed}, an 
absurdity results. It is
called by Gr\"{u}nbaum (following Weyl) ``Bernoulli's fallacy'': {\em 
He [Bernoulli]
treated the actually infinite set of natural numbers as having a {\em last} or
``$\infty$th'' term which can be ``reached'' in the manner in which 
an inductive cardinal
can be reached by starting from zero.}\footnote{A. Gr\"{u}nbaum, {\em 
o.c.}, pp.
130-131.} The error imposed on Zeno has therefore two wings that 
match the horns of the
presumed dilemma; one concerning `infinitely small' and one 
concerning `infinitely
large'.\footnote{Compare the sections ``The Deduction of Nullity of 
Size'' and ``The
Deduction of Infinity of Size'' in Vlastos's discussion of the problem, VLAS,
``Plurality'', p. 227 sq. and p. 233 sq.} The logical explicitation 
of the first wing has
been summarised nicely by Vlastos: {\em Since ``infinitely'' = ``endlesly'' and
``completion'' = ``ending'', it follows that ``the completion of the 
infinite division of
$x$ is logically possible'' = ``the ending of the endless division of
$x$ is logically possible.''} \/Its mathematical variant is the 
Bernoullian fallacy
inhaerent, according to Gr\"{u}nbaum, in e.g. Lee's and Tannery's 
construal of Zeno's
plurality-argument: {\em it is always committed when the attempt is 
made to use {\em
infinite divisibility} of positive [i.e., finite] intervals as a 
basis for deducing
Zeno's metrical paradox and for then denying that a positive interval can be an
infinitely divisible extension.}\footnote{A. Gr\"{u}nbaum, {\em 
o.c.}, pp. 131-132.} Of
course it is true that an infinite number of physically extended 
things lumped together
would be infinitely big! It is not even necessary to suppose that 
they be equal, as
Huggett\footnote{N. Huggett, {\em o.c.}, section 2.2.} does. For 
evidently I can build
unequal parts out of equal ones, as long as they stand in rational 
proportions to each
other, which is exactly what the decreasing sequence of cuts of the 
received view brings
about. So whether one construes Zeno's procedure as a directed division
or as one potentially throughout (as described above) will not even 
make a difference.
Once the sequence of  partitions is interrupted somewhere --- 
remember the meaning of
{\em a-peiron} --- and thus remains incompleted even after an 
infinite number of steps,
it will generate an infinity of parts with finite magnitude (unequal 
in the case of
directed division, equal in case the division was `throughout'). But 
Zeno denies exactly
this possibility of getting a countably infinite number of extended 
parts from a finitely
extended thing: {\em if they are many [things], they by necessity are 
as many as they
are, not more nor less. But if they are as many as they are, they 
will be bounded} \;
[{\small DK} {\scriptsize{29B 3}}]. In the interrupted case, the number will be
`bounded', i.e., finite, and the magnitude will be so too! Abraham 
calls this Zeno's
principle of {\bf the equivalence of the parts and the whole}, and as 
far as I can see it
is the only way to understand it commensurable with the content of the other
fragments.\footnote{W.E. Abraham, ``Plurality'', pp.
40-52. Contrary to KRS, who construe the second part of [{\small DK}
{\scriptsize{29B 3}}] as an objection to the first.}  A valid 
interpretation of `directed
division' with respect to this is hardly conceivable. Finally, 
following Abraham's
linguistic argumentation\footnote{Not to mention the fact that the 
direction imposed upon
the procedure in all likelihood is merely an artefact of our 
direction of reading! W.E.
Abraham, ``Plurality'', p. 42.}, this construal of Zeno's procedure 
can definitively be
ruled out, because it is perfectly possible to translate the 
substantivated participle in
the mentioned clause in [{\small DK} {\scriptsize{29B 1}}] as ``{\em 
any} projecting
[part]'', since ancient Greek does not discriminate between a 
definite and an indefinite
particle.\footnote{I agree with Abraham that the parallel with
$\pi\csup{\alpha}{-1}\varsigma \ \spir{o}{-1} \
\beta o\upsilon\lambda\acca{o}{-1}\mu\epsilon\nu o\varsigma$\; is 
relevant. Of the
individualising vs. the generalising --- ``to make a certain person 
or thing into the
representative of the whole species'' --- function of the article, 
very explicit examples
can be given. The translated quote stems from P.V. Sormani and H.M. 
Braaksma, {\em
Kaegi's Griekse Grammatica}, Noordhoff, Groningen, 1949, pp. 
112-113.} When one has no
content-loaden a priori in mind, `any' evidently fits in: it would 
simply mean any of the
symmetric parts obtained by a dichotomic division, irrespective of 
the level the
procedure attained. This is confirmed by the use of the verb
$\spil{\alpha}{-1}\pi\acca{\epsilon}{-1}\chi\epsilon\iota\nu$ in the 
previous line, which
conveys the idea `to be kept away from each other, to be separated 
actively' --- and
therefore to remain in contact, to be contiguous, like when one 
pushes his way through a
crowd --- rather than merely `to be at a distance'. This is why 
Vlastos in his paraphrase
of Zeno's arguments speaks of ``nonoverlapping 
parts''.\footnote{VLAS, p. 225.} The nice
and certainly not arbitrary antisymmetrical description of the 
totality of the procedure
when one takes [{\small DK} {\scriptsize{29B 1}}]: {\em for none will 
be the last
[eschaton] nor the one [heteron] will be unrelated to another one 
[pros heteron]}, and
contrasts it with [{\small DK} {\scriptsize{29B 3}}]: {\em For there 
will always [aei] be
others [hetera] in between [metaxu] of the beings, and there again 
others in between}, is
a case in point. The second wing also comes in two variants, though 
the difference is
between mathematics and physics, rather than between logic and mathematics. The
mathematical one is the familiar objection that the sum of an 
infinite series can be
finite, in case the series converges. The other variant goes back to 
Aristotle [e.g. {\em
Phys.} 263a(4-6)] and points in one way or another out that the 
execution of an infinity
of discrete physical acts within a finite stretch of time must be 
impossible. This
objection is at the origin of the recent discussion on 
``supertasks''.\footnote{For an
overview: J. P. Laraudogoitia,``Supertasks'', {\em STF},
http://plato.stanford.edu/archives/win2001/} These discussions can be 
relevant with
regard to Zeno, but only when one gives heed to certain precautions. 
I however want to
contest the  standard interpretation, and therewith reject the claim 
that Zeno committed
any such fallacy. It is important to stress here once again the {\em 
physical nature} of
Zeno's thought-experiment and the deictic realism inhaerent in it. It 
is as well
important not to introduce any presuppositions on behalf of Zeno, 
especially not those
that only came into existence in an attempt to deal with problems 
raised by him. To this
latter kind belong all mathematical and physical assumptions 
concerning the `continuous'
vs. the `discrete' nature of matter, space and time. {\em Zeno 
nowhere mentions space nor
time simply because these notions do not yet exist!}
\/Therefore he {\em  cannot} set up a dilemma on these praemisses. 
But then couldn't one
construct the dilemma on the basis of an utterance present in the 
Zenonian text? There is
no philological source-material upon which it can be based. Worse, of 
the necessarily
disjunctive structure underlying such an argument there is no trace. 
The contradiction
that should arise out of it is laid out already in the supposed 
praemisses, for Zeno uses
towards the end of [{\small DK} {\scriptsize{29B 1}}] twice the word
$\kappa\alpha\accg{\iota}{-1}$ {\em [and]} instead of the required
$\spilacca{\eta}{0}; \ \spilacca{\eta}{-1}\tau o\iota$ {\em [or]}, 
thus complying once
again to Simplicius's dictum that he intended nothing else but to 
show {\em that a man
who says that there is a plurality is stating something contradictory}.\\

Now, in order to throw more light on the nature of the problem at 
hand, it is crucial to
realise that there is a difference between the joining of two 
physical parts and the
addition of two mathematical line segments. The difference touches 
upon the far from
trivial problem of the relation between the mathematical and the 
natural sciences. In
Plato's {\em Phaedo} there is an extremely illuminating discussion of 
this relationship,
with reference to the ``contradictory'' results generated by the 
natural philosophy gone
before. Socrates complains that, since it allows for akin things to 
have contradictory
causes, while phaenomena clearly distinct become causally 
undistinguishable, its results
cannot be considered valid [{\em Phaedo}, 
100(e)-101(a,b)].\footnote{All references to this
dialogue are to the text in the LOEB$_3$-edition.}  He gives the example of the
difference between ``being two things'' and ``being a pair of 
things''; the latter a
formal, the former a physical fact. He also focuses on the relation 
between parts and
wholes in number theory, by comparing the generation of `two' out of 
`two ones' by by
adding them, bringing them together, and the separation of `two' into 
`two ones' by
dividing it [{\em Phaedo}, 97(a,e)]. The Zenonian influence is plain,
although scarcely discussed in the literature. There is an even more 
striking parallel
with Plato's {\em diairesis} -- the dissection of a concept into its 
constituting
contraries, the technique shoring up properly practiced dialectics, 
as applied in the
{\em Sophist}, the {\em Statesman}, and the {\em 
Philebus}.\footnote{The texts I
consulted are those in LOEB$_4$.} We will return in more detail to 
the {\em Philebus},
where the link between conceptual and mathematical {\em diairesis} is 
established, when
discussing our proposal for a faithful mathematical representation of 
Zeno's divisional
procedure. As is well known. Plato in the {\em Sophist} defines 
`dialectics' as the art
of making the proper distinctions between the forms that instantiate 
themselves in and
through particular things [{\em Sophist}, 253d(1-3)]. In that 
dialogue, moreover, the
difference established between contraries and contradictions --- 
between praedicative and
existential paradoxes --- lays out its ontological preconditions.\footnote{{\em
Platonists who doubt that they are spectators of Being must settle 
for the knowledge that
they are investigators of the verb `to be'.} G.E.L. Owen, ``Plato on 
Non-being'', {\em
Plato: a Collection of Critical Essays}, vol. i, G. Vlastos ed., 
Anchor/Doubleday, N.Y.,
1971, p. 223.} Now, harking back to something used already in the 
{\em Parmenides}, viz.
the dichotomic way of reasoning, the method by which {\em diairesis} 
should be applied is
demonstrated in the {\em Philebus} and the {\em Statesman} by means 
of examples. In the
latter dialogue, while trying to define the good statesman, Socrates 
and his friends find
out that the most long and cumbersome, but nevertheless the best way 
[{\em Statesman},
265(a)] to discover the specific forms instantiated in a thing is by 
systematically
dividing its concept in opposing halves, like `living/non-living' 
[{\em id.}, 261(b)],
`feathered/unfeathered' [{\em id.}, 266(e)] or `odd/even', instead of 
arbitrarily
separating off a part --- `Greeks' vs. `barbarians'; `ten' vs. `all 
other numbers' [{\em
id.}, 262(d-e)], say. This process ends when one bumps on 
undetermined parts, the {\em
stoicheia} or elements, that are not themselves capable of being 
specified further into
underlying parts [{\em id.}, 263(b)]. There are --- in the vocabulary 
of later times ---
no {\em differentia specifica} involved anymore, their plurality, if 
any, is  merely
numerical. The number of steps needed to reach from the original 
unity to this fully
determined level --- the proportion between part and whole --- then 
defines somehow the
original concept [{\em Philebus}, 16(d)]. This, however, is often not possible,
especially not when the praedicates are relative properties like `warm/cold';
`short/tall' \&c. Their opposites will run apart into their proper 
infinities and thus
into conceptual absurdity unless a limit, a boundary [{\em peras}] is 
imposed on them, in
order to find the good {\em measure} that guarantees their 
non-destructive, bounded
aequilibrium from the reality of the thing they describe sprang [{\em 
Statesman},
283(d,e)]. Aristotle states it in words that could not be more clear: 
{\em Plato, for his
part, recognises {\bf \em two infinities, the Large and the Small}
\/$\tau\accg{o}{-1} \ \mu\acca{\epsilon}{-1}\gamma\alpha \ 
\kappa\alpha\accg{\iota}{-1}
\ \tau\accg{o}{-1} \ \mu\iota\kappa\rho\acca{o}{-1}\nu$ [to mega kai 
to mikron]} [{\em
Phys.} 203a 15]. They can be discriminated indubitably, for with each 
goes a different
mode of realisation: {\em everything is infinite, either through {\bf 
\em addition} [i.e.
stepwise], either through {\bf \em division} [i.e. simultaneously] 
(...)} [{\em Phys.}
204a 6].\footnote{See J. Stenzel, {\em Zahl und Gestalt Bei Platon 
und Aristotles},
Teubner, Leipzig, 1933, p. 30 sq.; p. 60 sq.}\\

Aristotle himself discusses the problem directly in relation to Zeno, 
and puts it in a
way that clears the road for his notorious solution: the introduction 
of the difference
between {\em potential} and {\em actual} infinity: {\em For in two 
ways it can be said
that a distance or a period or any other continuum is infinite 
[apeiron], viz., with
respect to the partitions [diairesin] or with respect to the 
projecting parts [tois
eschatois]} [{\em Phys.} Z, 2, 263a (24-26)]. The vast majority of 
authors admits that
Aristotle did indeed attempt to solve PM (though not PP!) that way, 
although his solution
is generally dismissed as false or no longer mathematically relevant: 
{\em We won't
pursue this position, for the actual/potential distinction is not 
applicable to modern
mathematics.}\footnote{N. Huggett, {o.c.}, p. 40.} This is hardly 
acceptable a priori when
one realises that it throws us back to a problem already posed in 
Plato's {\em Phaedo}:
what does it mean to divide a continuous one into two parts? 
Aristotle comments: {\em For
whoever divides the continuum into two halves thereby confers a 
double function on the
point of division, for he makes it both a beginning and an end} [{\em 
Phys.} $\Theta$
\, 8, 263a (23-25)].\footnote{translations are after Wicksteed and 
Cornford in LOEB$_2$.}
And the points of partition serve as boundaries to the parts, be they 
potential or actual
in number, i.e., discernable and countable or not.\footnote{This 
relationship is exposed
admirably clear by Mary Tyles, {\em The Philosophy of Set Theory. An Historical
Introduction into Cantor's Paradise}, Dover, N.Y. 2004 [$1^{th}$ ed. 
1989], especially
pp. 10-31. The book discusses the foundations of set theory with the 
clash between
finitists and non-finitists on the foundations of mathematics as its 
vantage point.
Interestingly enough it has a logical counterpart as well. C. Vidal 
showed in his {\em
Ma\^{i}trise de philosophie:} ``Georg Cantor et la d\'{e}couverte des 
infinis'', Paris
I-Sorbonne, 2002-2003, p. 40, that a far from innocent inversion of 
logical quantifiers
is involved.} Zeno's Gedanken-experiment, by investigating the 
ultimate consequences of
the infinite divisibility of a physical being, hits on the said 
difference between the
physical and the mathematical and thus creates the preconditions 
necessary to establish
its various later forms. Poincar\'{e} discussed this crucial 
difference in a way that
will prove of great relevance to the mathematical representation we 
are going to build in
the remainder of this paper.\footnote{H. Poincar\'{e}, ``Le continu 
math\'{e}matique'',
{\em Reveu de m\'{e}taphysique et de morale}, {\bf 1}, 1893, pp. 
26-34. A reprint
entitled ``La grandeur math\'{e}matique et l'exp\'{e}rience'' in {\em 
SHYP}, pp. 47-60. He
deepens his ideas in ``La notion d'espace'', reprinted in {\em SVAL}, 
pp. 55-76.} Tannery
refers in this context to the critique of Protagoras on the geometers 
as quoted in
Aristotle [{\em Met.} B, ii, 997b(33)-998a(5)] : {\em (...) for as 
sensible lines are not
like those of which the geometrician speaks (since there is nothing 
sensible which is
straight or curved in that sense; the circle touches the ruler not at 
a [single] point, as
Protagoras used to say refuting the geometricians)}.\footnote{In a 
rare display of
negligence, DK destroy exactly the point of the argument by translating
$\kappa\alpha\nu\acca{o}{-1}\nu o\varsigma$ by {\em Tangente} instead 
of {\em ruler}! And
although Tannery places this discussion explicitly ``dans le cadre de 
Z\'{e}non'', he
concludes: {\em (...) la question est donc d'un autre ordre que 
celles soulev\'{e}es par
Z\'{e}non, et sa dialectique \'{e}tait impuissante \`{a} la 
r\'{e}soudre (...)}, the
reason being that also he takes Zeno's argument as constituting the 
dilemma exposed
above. It will be clear that I do not approve of either position: 
when understood
properly, Zeno does not have this problem at all. See [{\small DK} 
{\scriptsize{80B 7}}],
transl. in vol. II, p. 266; P. Tannery, {\em o.c.}, pp. 396-397. I 
used the translation
given by H. Tredennick, LOEB$_1$, Aristotle, {\em Metaphysics}, p. 
115. The [ ] are
mine.}  So, again, we cannot find salvation with logic, nor with space and
time\footnote{As KRS, p. 273, point out with regard to the Arrow Paradox.},
nor with the continuous and the discrete, just with Zeno's division 
throughout. The only
other thing we can rely upon is Parmenides's {\em Being-now}. This 
analysis imposes two
criteria that have to be fulfilled by whatever formal representation 
of Zeno's paradox:
{\bf constructivity} and {\bf simultaneity}. It indeed endorses a 
hitherto haedly
recognised line of approach. It is Abraham who (overlooking 
Parmenides, but reference to
him only strengthens the reasoning) opens up in a brilliant argumentation this
possibility, which comes down to accepting fully the consequences of 
the fact that Zeno's
division is a timeless division. To start with, let us follow what he 
has to say:
\begin{quote} {\em The objection that Zeno assumes the completion of 
an infinite task
assumes that, when he postulates that the being is divided {\bf \em 
through and through}
and so on infinitely, he is introducing end-products which logically 
cannot be further
divided, or he is assuming a least division beyond which there is no 
other, and so, a
last part. Now even though there is an actual infinite number of 
points in a line at any
of which points the line may be divided, the finite line does have 
terminal points
($\ldots$) But it would be a howler, committed by Johan Bernoulli, 
and decried by
Leibniz, that a terminal point would be ``the infinitieth point'' on 
the line. ($\ldots$)
To say that it is infinitely divided is no more than to say it 
actually has an infinite
number of points at {\em every} one of which it is divisible. The 
point to note is that
the infinite divisibility means not an infinite number of points of {\em
alternative}\footnote{i.e., stepwise. The italics are in the 
original. The bold further
down the quotation is mine.} division (such that the alternatives are 
inexhaustible) but
rather an infinite number of points of simultaneous division. The 
points of division,
being points on the being, belong to it not alternatively, but {\bf {\em
simultaneously.}} It is this simultaneity (and not a process) which 
is articulated by the
postulate of the complete division. It is clear that if Zeno's 
complete division thus is
a {\bf \em cardinal} completion rather than an {\bf \em ordinal} 
completion, the infinite
division of the given being does not imply a {\em last} division or 
{\em last} part, any
more than the simultaneity of the points on a line imply an infinitieth
point.}\footnote{W.E. Abraham, ``Plurality'', p. 48. My bold.}\end{quote}

\noindent The central point thus is that by being simultaneous (i.e., 
by occuring in no
lapse of time) the completion by division is at once, cardinal, not 
stepwise, ordinal.
This refers to the aritmetic of ``infinity'', developed towards the 
end of the nineteenth
century by G. Cantor, who proved that different kinds of infinity 
exist which he called
{\em cardinal} and {\em ordinal} infinity.\footnote{G. Cantor, 
``Beitr\"{a}ge zur
Begr\"{u}ndung der transfiniten Mengenlehre'', {\em Gesammelte 
Abhandlungen} [{\em GA} in
what follows], E. Zermelo Ed., Georg Olms Verlag, Hildesheim, 
1932/1962, pp. 312-356.}
Cardinality (``Zahl") expresses the total number of elements of a 
given set; ordinality
(``Anzahl") concerns the way these elements can be ordered stepwise. 
He furthermore
showed that the set of natural numers $\nn$ and the set of finite 
fractions $\mathbb{Q}$
possess the same {\em countably infinite} cardinality, which he 
labelled Aleph-null
($\aleph_0$). The basic ordinality that goes with it he called 
$\omega$. It is possible
to construct an infinite sequence of $\omega$'s which comprise 
different ordinalities that
all belong to the same {\em number class}, i.e., discern levels of 
equal cardinality:
$\omega_i = Z(\aleph_1)$; $\omega'_i = Z(\aleph_2)$, \&c. The 
difference only becomes
relevant from the moment the numbers implied are infinite. Consider 
for instance theset
of natural numbers ordered in the traditional manner $\{0, 1, 
2,\ldots \}$ and ordered
alternatively $\{1, 2,\ldots , 0\}$. The cardinal number of these 
sets will be equal;
their ordinal numbers will be respectively $\omega$ and $\omega + 1$. 
Cantor also showed
by means of his famous diagonalisation argument that
$2^{\pw{\nn}}$, the cardinality of the set of real numbers
$\rr$\footnote{F. Hausdorff, {\em Mengenlehre}, Chelsea Publishing 
Company, N.Y.,
1949, pp. 62-64.} --- the mathematical face of the idea of continuity 
---, is an infinity
{\em bigger} than the cardinality of $\nn$ or $\mathbb{Q}$: they do 
{\em not} \/belong to
the same number class.\footnote{G. Cantor, ``\"{U}ber eine elemetare Frage der
Mannigfaltigkeitslehre'', {\em GA}, pp. 279-281.} The infinity of the 
continuum is {\em
uncountable}. Abraham's intention with respect to Zeno's division now 
becomes more clear.
He apparently claims that Zeno's procedure {\em at once} generates 
the {\em uncountable}
cardinality of the continuum instead of the merely {\em countable} 
one which would be
attainable when one interpretes Zeno's procedure {\em stepwise}. This 
reminds us of the
viewpoint also defended by authors like Tannery and Luria, castigated 
by Gr\"{u}nbaum for
reasons given in the following comment: \/{\em it is essential to 
realize that the
cardinality of an interval is not a function of the length of that 
interval. The
independence of cardinality and length becomes demonstrable by 
combining our definition
of length with Cantor's proof of the equivalence of the set of all 
real points between 0
and 1 with the set of all real points between {\em any} two fixed 
points on the number
axis.}\footnote{A. Gr\"{u}nbaum, {\em o.c.}, p. 127.} And it is 
exactly this error that
Gr\"{un}baum predicates on Zeno's plurality paradoxes: according to 
him, Zeno maintains
that a longer part ``contains more points''. This is nothing but 
another variant of the
`Bernoullian fallacy'' which he believes also to be present in the 
paradoxes of motion.
Now it is precisely Abraham's aim to show that this error is {\em 
not} present in Zeno's
plurality-argument, for it would imply a divisional process {\em 
through time,}\/ and
Zeno nowhere mentions `time'. But this does not yet suffice to 
completely solve the
riddle, because it leaves untouched a related problem raised by 
Gr\"{u}nbaum: how to
construct an extended object out of parts with no extension 
whatsoever? Zeno cannot just
intend that his division generates the unextended points composing 
the continuous line,
for then only the `being nothing'-part of his paradox would remain. 
Zeno's divisional
procedure generates an infinity of partitions and a {\em different} 
infinity of parts at
once and jointly {\em in any part independent of its length}.  The 
simultaneous and
uninterrupted [{\em apeiron}] repetition of partitions is clearly 
stated by Zeno himself:
  {\em to say this once is as good as saying it forever [aei]}. The 
result is an infinity
of parts which must at least be dense: {\em for none will be the last 
nor one will be
apart from another one} and {\em For there will always be others in 
between of the
beings, and there again others in between}. We know moreover that it 
results in parts
with and without magnitude, whence we may safely assume that this 
strange fact is somehow
connected to the presence of different infinities: we do not only 
have density but also
arbitrary length! It is tempting to conclude that Zeno's way of 
putting his problem
betrays an intuitive awareness of the notion of infinite cardinality, 
which moreover would
render Aristotle's bewildering step to discriminate between two 
infinities, potential and
actual, not only comprehensible, but places it in a glaring light. 
This is confirmed by
Aristotle's own interpretation: {\em For in two ways it can be said 
that a distance or a
period or any other continuum is infinite [apeiron], viz., with 
respect to the partitions
[diairesin] or with respect to the projecting parts [tois 
eschatois]}, although Aristotle
in all likelihood would deny that Zeno himself realised this. But a 
drawback of our
interpretation up to now is that it does not fit into the received 
view on the paradoxes
of motion, of which some are incontestably based on dichotomy as 
well, as again the
Stagirite diagnoses: {\em This argument [the Achilles] is the same as 
the former which
depends on dichotomia} [{\em Phys.} Z, 9, 239b (20-21)]. The received 
view on Zenonian
plurality on the contrary does, at least partially, comply to this 
criterion. I can
circumvent this objection only by devising an interpretation that 
gives us PM as the
result of a division through and through. This is precisely the 
nature of my claim with
respect to them, as I hope to convincingly argue towards the end of 
this article.

\bigskip

{\sc A simple mathematical representation}.--- In accordance with 
Zeno's already
mentioned realism towards everyday phaenomena like plurality and 
change, we shall develop
a mathematical representation of Zeno's divisional procedure
$\mathfrak{Z}$ that is as simple as possible, that does not make any 
presupposition other
than those explicitly retracable to our philosopher, and is 
applicable to all his
paradoxes alike. So, then, what do we know about Zeno's division?

\begin{tabbing}
	\= i) it goes on ``infinitely''; \\
	\> ii) it is symmetrical (proportion 1:2); \\
	\> iii) it results somehow in parts with and without magnitude; \\
	\> iv) it is a simultaneous procedure; \\
	\> v) it pops up in all his paradoxes.
\end{tabbing}

\noindent  What we need is a procedure whereby ``two different kinds 
of infinity are
generated simultaneously'', and in which ``the points of division are 
used twice''. Let
us take as a model object a measuring rod
$M$, of which we conventionally call the left-hand side 0 and the 
right-hand side 1, and
which we arbitrarily set equal to unity.\footnote{I owe this 
illustrative example to an
electronic discussion with J. Helfand.} This we can do, given that 
the specific length of
the object by no means influences the argument. What will happen when 
we Zenonian-wise
divide this model-object through and through? We can see this 
division as a kind of cell
division, whereby in each generation the number of parts generated 
doubles: when
generation $n$ has
$2^n$ parts, generation $n + 1$ will have
$2^{n+1}$ parts. And we know that $n \pijltje \infty$. This can be represented
graphically by a divisional tree. The cells function at every step in
two different ways: as partitions in generation $n-1$ and as parts in 
generation $n$. In
every generation, generation by addition and generation by division 
thus coincide. The
number of partitions should at least reach
$n \leq \omega$, with $\omega$ the basic ordertype of $\nn$:\\

\begin{center}
\includegraphics[width=.6\textwidth,angle=0]{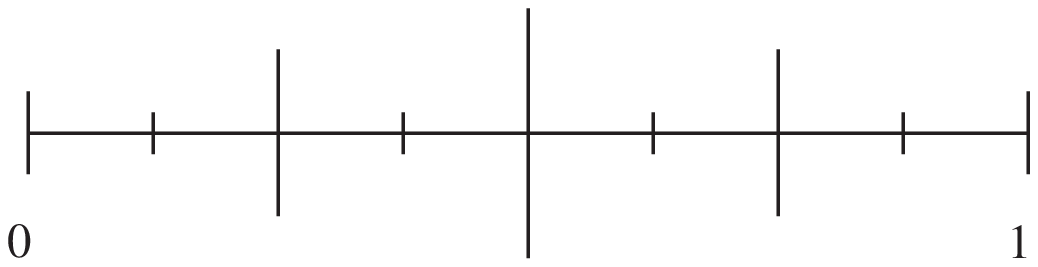}
\end{center}

\vspace{-8mm}

\begin{center}
\includegraphics[width=.6\textwidth,angle=0]{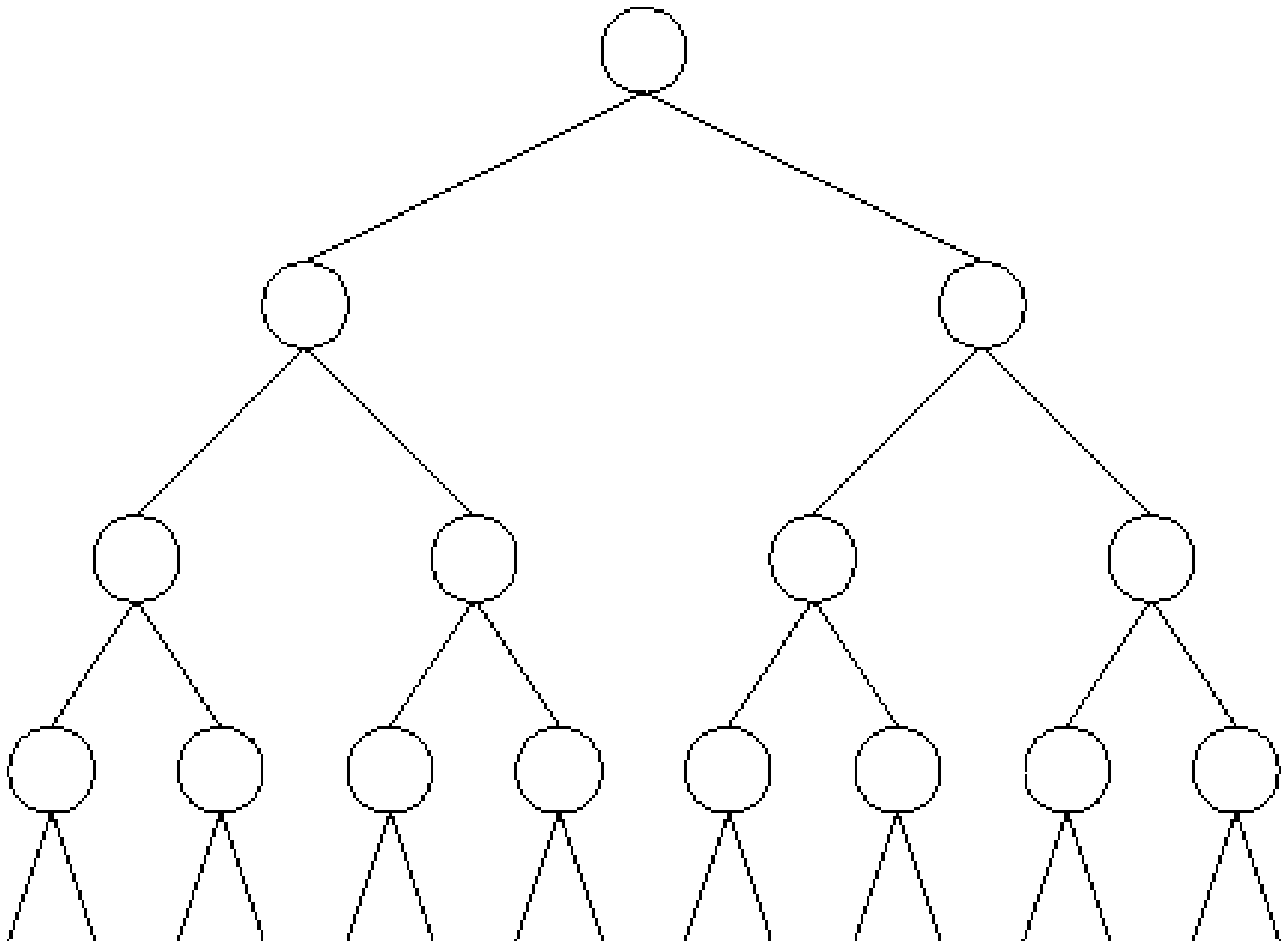} \\ 
\vspace{-6mm} {The
measuring rod
$M$ and the tree of cell-divisions.}
\end{center}

\noindent The simultaneous division-lattice stands for the totality 
of $\mathfrak{Z}$, and
thence equally bears a twofold interpretation: by looking at the 
``last row" at generation
$\omega$ and {\em dividing once more}, or by mapping the whole procedure on the
$\nn$-lattice by relating the number of partitions to number of 
cells, and the number of
parts to the number of paths between them. The first procedure comes 
intuitively closer
to what Zeno intended, but inevitably introduces an inappropriate 
notion of time, a flaw
which is absent from the second one. Anyhow, the result in both cases 
is the same: the
number of partitions $\mid\!\!\nn\!\!\mid$ with its concommittant number of
parts $\mid\!\!\rr\!\!\mid = 2^{\footnotesize{\nn}}$ are generated 
simultaneously and with
a clearly {\em different} status. So I am {\em not} saying that the 
result of Zeno's
procedure is merely a set with the cardinality of the continuum, like 
the Euclidean
straight line. That would, moreover, be nothing new, as I indicated 
when discussing
Tannery's and Abraham's views. It nevertheless is already remarkable 
enough that such a
construction of the continuum should be possible at all. Attempts to 
do something akin to
it have been proposed in the context of constructive or intuitionistic
mathematics.\footnote{My attention was drawn to intuitionistic mathematics by
the discussion in M. Tiles's book, {\em o.c.}, p. 90-94.} In order to make the
difference with our own procedure completely clear, we will devote 
ourselves now to a
short and non-rigorous discussion of the approach on the basis of 
Brouwer's dissertation
and Troelstra's fundamental work.\footnote{A.S. Troelstra, {\em Principles of
Intuitionism}, Springer, Berlin/Heidelberg, 1969, p. 22, 52; pp. 
57-64.} Briefly stated,
{\bf intuitionistic mathematics} was conceived by L.E.J. Brouwer in a 
reaction to the
discovery of the logical and set theoretic paradoxes at the beginning 
of the twentieth
century. The basic philosophical ideas shoring up his d\'{e}marche 
are 1) mathematics
{\em precedes} logic, not the other way around; and 2) classical logic fails in
universally applying the principle of the excluded middle (TND), 
which states that {\em
all} logically valid propositions have either the values ``true'' or 
``false'' attached
to them. Both principles are intimately connected to Brouwer's 
attitude towards the
mathematical infinite. To talk about ``all'' or ``for every'' when 
the number of objects
considered is infinite amounts into absurdity according to him. He 
holds that in
mathematics the infinite can only be {\em potential}, i.e., so as to 
be stepwise
approached but never reached, while the {\em actual} infinite of the 
standard view of the
continuum cannot have real existence. In his {\em Proefschrift} 
[dissertation], Brouwer
writes: {\em Het continuum als geheel was ons echter intu\"{\i}tief 
gegeven; een opbouw
ervan, een handeling die ``alle'' punten ervan ge\"{\i}ndividualiseerd door de
mathematische intu\"{\i}tie zou scheppen, is ondenkbaar en 
onmogelijk. De mathematische
intu\"{\i}tie is niet in staat anders dan aftelbare hoeveelheden 
ge\"{\i}ndividualiseerd
te scheppen, n.l. z\'{o}, dat men een proc\'{e}d\'{e} geeft, dat elk 
element der
verzameling na een eindig aantal operaties genereert.}\footnote{{\em 
The continuum as a
whole was given to us intuitively however; to construct it by means 
of an act of
mathematical intuition that would create ``all'' of its points individually, is
unthinkable and impossible. Mathematical intuition can only create 
countable quantities
individually, i.e., such, that a procedure is given which generates 
every element of the
collection after a finite number of operations.} L.E.J. Brouwer, 
``Over de grondslagen
der wiskunde. Academisch Proefschrift [1907]'', in D. Van Dalen, {\em 
L.E.J. Brouwer en
de grondslagen der wiskunde}, Epsilon, Utrecht, 2001, p. 73. [My 
translation.]} The
question then becomes how much standard mathematics  we can get back 
on such a constrained
conceptual basis. Brouwer himself advanced the idea that the continuum could be
constructed by means of finite sequences --- a countable number of 
``duaalbreuken,
ternaalbreuken'' [dual, ternal fractions, i.e., rationals written in 
binary or ternary
form] --- which approach points
$P$, the ``dubbelpunten'' [double points]\footnote{ L.E.J. Brouwer, 
{\em o.c.}, pp. 56.}
demarcated by neighbourhoods which remain outside the ``schaal'' 
[scale] determined by the
``voortschrijdingswet'' [Law of Progression] of the 
sequence.\footnote{L.E.J. Brouwer,
{\em o.c.}, pp. 43-45.} It is as if one covers the given continuum 
with another one full
of tiny holes [``lacunes''] which one can then ``samendenken'' [think together]
by the process of approximation.\footnote{{\em Zo kan een gegeven 
continu\"{u}m door een
ander continu\"{u}m met lacunes worden overdekt [covered]; we 
behoeven daartoe op het
eerste continu\"{u}m maar een ordetype $\eta$ te bouwen, dat het niet 
overal dicht [dense]
bedekt en vervolgens bij dat ordetype $\eta$ het continu\"{u}m te 
construeren; we kunnen
dan altijd een punt van het tweede continu\"{u}m identiek noemen met 
het grenspunt van
zijn benaderingsreeks [the limit-point of its approaching sequence] 
op het eerste
continu\"{u}m.} {L.E.J. Brouwer, {\em o.c.}, p. 73.} The ``ordertype 
$\eta$'' belongs to
Cantor's {\em first number class}, as Brouwer makes clear on pp. 42-43
[cfr. ft. 72 infra].} This process can be represented by a ``boom'' 
[tree]\footnote{L.E.J.
Brouwer, {\em o.c.}, pp. 73-76.}, to the decreasing parts of which 
one can associate a
binary representation.\footnote{H. Weyl, {\em Philosophie der Mathematik und
Naturwissenschaft}, M\"{u}nchen, 1926, p. 43; p. 51.} Troelstra 
formalises these ideas;
the basic notions required are a {\em species}, a {\em sequence} and 
a {\em tree} or {\em
spread}. Species are the intuitionistic analog for classical sets. 
They comply to a
specific form of the {\em axiom of comprehension}, predicated on a 
precise interpretation
of the concept of {\em well-definedness} based on the idea of 
provability. (Again, we
shall not forget that intuinionism originated in the wake of the 
discovery of the famous
set theoretic paradoxes.) Now let a sequence be a mapping, {\em i.e. 
a {\bf \em process}
which associates with every natural number a mathematical object 
belonging to a certain
species
$X$}: $(\nn)X: \, \nn \pijlr X$.\footnote{A. Troelstra, {\em o.c.}, 
p. 16.} These
sequences are introduced in order to develop an intuitionistic theory 
of real numbers.
Given a map $\xi$, the species $X$ can be {\bf represented} by a spread with an
aequivalence relation
$\sim$ if $<\!a, \, \xi\!>$ is a spread, and if
$X^{\ast}$, the species of aequivalence classes with respect to 
$\sim$ can be mapped
bi-uniquely (injectively) onto $X$. So $X$ can be represented by 
$<\!a, \, \xi\!>$
  via its partition $X^{\ast}$. {\em The spread represents a set of 
nodes of a growing
``tree'' of finite sequences of natural numbers, with branches 
directed downwards; the
topmost node corresponds always with the empty sequence.}\footnote{A. 
Troelstra, {\em
o.c.}, p. 52.} The tree may be branched up to an arbitrarily high 
level. It will be clear
that this construction does {\em not} fulfill the criteria which 
emerged out of our
analysis of Zeno's procedure $\mathfrak{Z}$. For, as Troelstra 
himself unambiguously
states, his procedure implies a {\em process} consisting of a finite 
though unspecified
number of steps. And in a related context, Brouwer himself says that 
he sees the sequence
of divisions generated by it as taking place {\em door den tijd 
[through time]}. The
open-endedness of the procedure is of course a necessity if one wants 
to approximate the
real numbers by means of it.\footnote{H. Weyl, {\em o.c.}, p. 44.} 
But how then to avoid
the Bernoullian fallacy? By rehabilitating the Aristotelian notion of the {\em
potentially} infinite, the idea of an infinite process that 
theoretically could, but in
reality never will, be concluded. {\em The totality of all {\bf
\em potentially infinite} sequences of rationals (...) is a 
potentially infinite whole
containing members which are never fully determinate, for they 
themselves are potentially
infinite, and always {\bf
\em actually finite} but necessarily incomplete.} As a matter of 
fact, {\em since the
process of focusing is not, and can never be complete, the continuum 
is never resolved
into points, only into {\bf \em ever smaller regions}.}\footnote{M. 
Tiles, {\em o.c.},
pp. 91-92.} Hermann Weyl, who stresses exactly this aspect of 
Brouwer's approach, links
it to Plato: {\em Brouwer erblickt genau wie Plato in der 
Zwei-Einigkeit die Wurzel des
mathematischen Denkens.} Plato's
$\spil{\alpha}{-1}\acca{o}{-1}\rho\iota\sigma\tau o\varsigma \
\delta\upsilon\acca{\alpha}{-1}\varsigma$ {\em [indeterminate 
twoness]} is the principle
of generation underlying {\em diairesis [conceptual division] } by 
means of which both
the Large and the Small come about. As Aristotle explains in Book
$N$ of the {\em Metaphysics}, it is closely related to Platonic number
theory.\footnote{This fact is in itself quite undeniable and its 
recognition does not
imply a position in the controversy surrounding the ``T\"{u}binger Schule''
interpretation of Plato's philosophy. Traces of
Platonic {\em diairesis} can be found back in the earliest dialogues, 
and arguably
contributed to the development of Plato's theory of Forms, at least 
according to M.K.
Krizan, ``A Defense of {\em Diairesis} in Plato's {\em Gorgias}, 
463e5-466a3'', {\em
Philosophical Inquiry}, XII(1-2), 1-21. For a recent and moderate 
overview of the issues
at stake, see D. Pesce, {\em Il Platone di Tubinga, E due studi sullo 
Stoicismo},
Paideia, Brescia, 1990.} We will see later that Weyl's account in 
fact only reckons the
`Large'-part of Platonic {\em diairesis}. It is outside our scope to 
dwell on this more
elaborately here, but the parallel with Zeno's analysis will nevertheless be
clear: Plato's Large-and-Small
$\tau\accg{o}{-1} \ \mu\acca{\epsilon}{-1}\gamma\alpha \ 
\kappa\alpha\accg{\iota}{-1}
\ \tau\accg{o}{-1} \ \mu\iota\kappa\rho\acca{o}{-1}\nu$ 
\footnote{Cited by Aristotle in
e.g. [{\em Met. A}, 987b(20)]. Cfr. J. Stenzel, {\em Zahl und Gestalt 
Bei Platon und
Aristotles}, Teubner, Leipzig, 1933 (2nd ed.), p. 6.} and Zeno's 
large[s]-and-small[s]
$\mu\acca{\epsilon}{-1}\gamma\alpha\lambda\alpha \ 
\kappa\alpha\accg{\iota}{-1} \
\mu\iota\kappa\rho\accg{\alpha}{-1}$ both refer to the
paradoxical result of an infinite and through and through 
division.\footnote{Although I
found it nowhere discussed from this perspective in the literature. 
The singular in
Plato's expression stems from the fact that he applies the 
paradoxical property to the
abstract principle by which division is obtained, while with Zeno the 
plural directly
refers to its result.}  Plato's metaphysical preoccupation with 
non-paradoxical plurality
however forces him to enlarge the `here'-part of Zeno's deictic 
standpoint into an in
principle unlimited number of partial perspectives that cover the 
totality of the
divisional tree. Thanks to an idea of M. Serfati, it will be possible 
to formally
distinguish the set of perspectives ({\em l'ensemble des points de 
vue}) from the tree
itself.\footnote{M. Serfati, ``Quasi-ensembles d'ordre $r$ et approximations de
r\'{e}partitions ordonn\'{e}es'', {\em Math. Inf. Sci. hum.}, {\bf 
143}, 1998, pp. 5-26.}
For the human mind, however, the relocation of the infinite to the 
realm of the Ideas
makes that the part of the tree accessible to it always remains 
partial, while with Zeno
the totality of the divisional process is given in the here-and-now. 
But again with
Plato, being ideal does not mean being unreal, quite the contrary! 
That is why the
Platonic spectator of Being in its One/Many-marked appearance retains 
a definitely
paradoxical flavour.\\

\noindent {\sc Some preliminary notions}.-- In order to fully 
establish $\mathfrak{Z}$ we
need, I repeat, a procedure whereby ``two different kinds of infinity 
are generated
simultaneously'', and in which ``the points of division are used 
twice''. The two kinds
of infinity are {\em simultaneously} there and have a {\em different} 
status: as
partitions and as parts. Because of the constructive simultaneity, they do not
only come both about at once, but pop up as well without anything 
else intervening
between them. This latter point gives us a clue as to where to look 
for a viable
mathematical approach, for in set theory a very important 
hypothetical theorem deals
exactly with this type of situation. This is {\bf Cantor's continuum 
hypothesis}. But
before we can appropiately make use of it ourselves, we first have to 
turn our attention
to the question that underlies it. Let us briefly review its 
mathematical and conceptual
content. We saw that Zeno's plurality paradox can be seen as a 
divisional procedure for
constructing the continuum as instantiated in a finitely extended 
body representable as a
part of the real line. It would be nice to know how much `parts' are 
exactly needed to
build it up. This is precisely the zest of Cantor's continuum 
problem: {\em How many
points are there on a straight line in Euclidean space?}\footnote{K. 
G\"{o}del, ``What is
Cantor's continuum problem?'', {\em BPPM}, p. 470.} We already saw 
that multiple infinite
cardinalities exist. They are represented by indexed Hebrew 
characters called alephs
($\aleph_i$). The continuum hypothesis is the question: what is the cardinality
$\aleph_c$ (or $c$, in a to-day more common notation) of the 
continuum? Let us take a look
at it more closely. Cantor defined the magnitude of a set by means of its {\em
cardinality}. Two sets are aequivalent, that is, a bijective 
relationship between their
elements exists, when their cardinalities are equal.\footnote{{\em 
unter eine {\em
``Belegung von
$N$ mit $M$"} verstehen  wir ein Gesetz, durch welches mit jedem Elemente
$n$ von $N$ je ein bestimmtes Element von $M$ verbunden ist 
($\ldots$) Das mit $n$
verbundene Element von $M$ ist  gewisserma{\ss}en eine eindeutige 
Funktion von $n$ und
kann etwa mit
$f(n)$ bezeichnet  werden ($\ldots$)}; G. Cantor, ``Beitr\"{a}ge zur 
Begr\"{u}ndung der
transfiniten Mengenlehre'', {\em GA}, p. 287.} But there is an other 
number that
expresses a basic characteristic of a set, its {\em ordinality}. 
Cardinality expresses
the simultaneous totality, the number of elements present of a given 
set; ordinality
concerns the order by way of which these elements are generated stepwise. {\em
``M\"{a}chtigkeit'' oder ``Kardinalzahl'' nennen wir den 
Allgemeinbegriff, welcher mit
Hilfe unseres aktiven Denkvermogens  dadurch aus der Menge $M$ 
hervorgeht, da{\ss} von
der Beschaffenheit ihrer verschiedenen  Elemente {\em m} und von der 
Ordnung ihres
Gegebenseins abstrahiert wird.}\footnote{G. Cantor, ``Beitr\"{a}ge'', 
GA, p. 282.} In the
finite case there is no problem: a basket containing eight apples is a set with
cardinality eight, its ordinality being eight as well, given that we 
count the apples one
by one. We can add one apple; now we have nine. We can get to any 
arbitrarily big natural
number this way. This is how Cantor constructs {\em potential} 
infinity. He furthermore
showed that the set of natural numers $\nn$ and the set of finite 
fractions $\mathbb{Q}$
possess the same {\em countably infinite} cardinality, which means 
that, although their
elements can be enumerated, they can never be exhausted by means of 
enumeration. He
labelled it Aleph-null ($\aleph_0$). Thus the set of natural numbers 
$\nn$ and the set of
finite fractions
$\mathbb{Q}$ are aequivalent. The basic ordinality that goes with it 
he called $\omega$.
It is possible to construct an infinite sequence of different ordinal numbers
$\omega_i$ that all belong to that same {\em number class} of 
cardinality $\aleph_0$. The
difference between ordinals becomes relevant from the moment the 
numbers so defined are
infinite. Consider for instance the set of natural numers ordered in 
the traditional
manner $\{0, 1, 2,\ldots \}$ and ordered alternatively $\{1, 2,\ldots 
, 0\}$.  The
cardinal number of these sets will be equal; their ordinal numbers 
will be respectively
$\omega$ and $\omega + 1$. The set \,
$\{2, 4, 6,\ldots ; 1, 3, 5,\ldots  \}$ (even and odd numbers) has 
ordinality $2\omega$.
The number of possible rearrangements is clearly unlimited. 
Nevertheless, they all belong
to the number class of $\aleph_0$. Now it is possible to build an 
infinite sequence of
ordinalities by application of two different ``Erzeugungsprinzipien'' 
or Principles of
Generation\footnote{G. Cantor, ``Beitr\"{a}ge''', {\em GA}, pp. 
312-356. Accessible
treatments can be found in R. Rucker, {\em Infinity and the Mind. The 
Science and
Philosophy of the Infinite}, Princeton University Press, Princeton, 
1982 [1995], p. 65
sq, p. 223; and M. Tyles, {\em o.c.}, pp. 104-107, of whom I borrow the English
terminology.}:

\smallskip

\begin{itemize}
	\item First Principle of Generation: {\em Addition
							of a unit to 
the previously formed number};
	\item Second Principle of Generation: {\em Definition of the
							limit of 
every infinite sequence of numbers as the nearest limit value 
[supremum]
							outside of it}.
\end{itemize}

\smallskip

\noindent Thus a new number is created justly outside the infinite 
sequence of already
existing ordinalities with cardinality $\aleph_0$, which represents an infinity
necessarily {\em greater} than $\aleph_0$. Cantor calls this new, 
bigger infinite
cardinality naturally enough $\aleph_1$. The first infinite 
cardinality $\aleph_0$ defines
the second number class $Z(\aleph_1)$.\footnote{The first number 
class $Z(\aleph_0)$
comprises the finite ordinalities embodied by the set $\nn$.} The new 
basic ordinality
can now be used to construct a new sequence of ordinalities, all belonging to
$Z(\aleph_2)$, the thirth number class. A sequence of ever
greater alephs can be built by systematically applying the two 
principles of generation
so that no upper bound can be assigned to it; with each aleph a new, 
bigger cardinality
aequivalent to the smallest ordinality of its respective number 
class: $\omega_i =
Z(\aleph_{i+1})$; $\omega_{i+1} = Z(\aleph_{i+2})$, \&c. Cantor 
proves by means of an
argument known as {\em diagonalisation} that indeed $\aleph_0 <
\aleph_1$. Another way to prove it is Cantor's powerset
Theorem:  $E \prec {\cal P}(E).$\footnote{F. Hausdorff, {\em 
Mengenlehre}, pp. 56-57;
G.S. Boolos and R.C. Jeffrey, {\em Computability and Logic}, CUP, 
Cambridge, 1974/1989,
p. 1 sq.} From his powerset theorem Cantor infers the general theorem 
that for each and
every cardinal number, as well as set of cardinal numbers, an 
immediate successor exists:
$\aleph_n < \aleph_{n+1}$, with $n \in
\mathbb{N}$\,).\footnote{F. Hausdorff, {\em Mengenlehre}, pp. 67-68.} 
It will be clear
that one thus not only obtains ascending sequences of ordinal, but 
also of cardinal
numbers $\aleph_0, \, \aleph_1, \, \aleph_2
\ldots$!\footnote{Cantor exposes his method in a succinct and clear 
way in a letter to
Dedekind (1899). A translated version in J. Heyenoort, {\em From 
Frege to G\"{o}del},
Harvard, 1967, pp. 113-117.} The sequence of {\em all} possible 
ordinals Cantor called
$\Omega$; the sequence of all possible cardinals $\aleph_0, \, \aleph_1, \,
\aleph_2 \ldots$ he called $\beth$.\footnote{{\em Idem}, {\em FTG}, 
pp. 113-117.}
These objects are themselves not sets, for ``thinking them as a 
whole'' leads to
paradoxes.\footnote{In the case of
$\Omega$ the Burali-Forti paradox: if the set of ordinals is well 
ordered (i.e., every
segment has a least element), it has an ordinal, which is at the same 
time an element of
this set and greater than any of its elements.  It is a variant of 
the more familiar
Russell paradox. G. Cantor, ``Letter to Dedekind'', in: J. Van 
Heijenoort, {\em From
Frege to G\"{o}del} [FTG in what follows], pp. 115-117. See also the 
paper by C.
Burali-Forti, {\em id.}, pp. 104-112. This fact --- {\em diese 
omin\"{o}se ``Menge W''}
[E. Zermelo] --- was the source of a vivid controversy at the 
beginning of the last
century. In addition to the papers present (in translation) in {\em 
FTG}, one will
profitably consult the volume composed by G. Heinzmann, {\em 
Poincar\'{e}, Russell,
Zermelo et Peano. Textes de la discussion (1906-1912) sur les fondements des
math\'{e}matiques: des antinomies \`{a} la pr\'{e}dicativit\'{e}}, 
Blanchard, Paris, 1986,
where a number of relevant but sometimes less well known papers are 
collected in their
original form. The Zermelo-quote stems from that source, p. 119. A 
relation to the axiom
of choice is exposed in M. Potter, {\em Set Theory and its 
Philosophy}, Oxford University
Press, Oxford, 2004, p. 243 sq.}  Such collections, let us say to 
`big' to be thought as a
whole, Cantor termend {\em inconsistent multiplicities}. Every 
possible cardinal number
will be represented somewhere in the sequence $\beth$. In order to 
prove this latter
point, one needs the often contested, but critical, {\em Axiom of Choice}
[AC].\footnote{K. G\"{o}del, ``What is Cantor's continuum problem?'', 
{\em BPPM}, p.
471.} Intuitively, the AC allows that {\em a {\bf \em simultaneous} choice of
distinguished elements is in principle always possible for an arbitrary set of
sets.}\footnote{E. Zermelo, ``A new proof of the possibility of a 
well-ordering'', {\em
FTG}, p. 186.} Zermelo, in his edition of Cantors collected work, 
castigates him for
failing to have realised this, thus committing a variant of the 
familiar Bernoullian
fallacy: {\em Only through the use of the AC, which postulates the 
possibility of a
simultaneous choice (...)}\footnote{E. Zermelo in FTG, p. 117, ft. 3, 
to Cantor' letter
to Dedekind reprinted therein.} Although essential to its logical 
viability, this
simultaneous character of AC is often forgotten or dismissed, while 
Zermelo himself never
gets tired of stressing its importance.  We now possess the 
instruments to reformulate
our original question: which one of the alephs in our sequence is 
equal to the number of
points $2^{\aleph_{0}}$ on the continuous line? Where exactly in the 
totality of all
possible cardinalities this number is to be situated?\footnote{{\em 
Il faut distinguer
entre {\em l'hypoth\`{e}se du continu} et le {\em probl\`{e}me du continu
(Kontinuumproblem)}, qui consiste \`{a} d\'{e}terminer la place 
occup\'{e}e par le
continu parmi les {\em alephs}, c.\`{a}.d.
\`{a} d\'{e}terminer le nombre ordinal
$\alpha$ \, pour lequel \,
$2^{\aleph_0} = \aleph_{\alpha}$}, explains W. Sierpi\'{n}ski, {\em 
Hypoth\`{e}se du
Continu}, Z Subwencji Funduszu Kultury Narodowej, Warsawa/Lw\'{o}w, 
1934 (reprinted in
1956 by the Chelsea Publishing Company), p. 5.}  Cantor demonstrated 
by means of
diagonalisation that $2^{\aleph_0} > \aleph_0$, and advanced the 
hypothesis that it
should be its immediate successor,
$\aleph_1$. This is Cantor's famous Continuum Hypothesis 
(CH)\footnote{K. G\"{o}del,
``What is Cantor's continuum problem?'', {\em BPPM}, p. 472.}:


\[ 2^{\aleph_{0}} =
\aleph_{1} \]

\noindent Every infinite subset of the continuum has the power either 
of the set of
integers, or of the whole of the continuum, with no intermediate cardinalities
intervening. As could be expected, the hypothesis can again be 
generalised, so as to
comprise all possible cardinalities: the Generalised Continuum 
Hypothesis.\footnote{Maar:
$\mid\!\!2^{\aleph_{0}}\!\!\mid\;$ kan iets anders zijn dan 
$\aleph_1$!} We will
notwithstanding limit ourselves to CH, for the nature of our problem 
is demarcated by
Zeno clearly to the realm of the physical, even when encountered in 
its mathematical
expression: the reality of plurality and motion, as problematised in 
the clash between
continuous and discrete: {\em the objects of transfinite set theory 
(...) clearly do not
belong to the physical world and even their indirect connection with 
physical experience
is very loose.}\footnote{K. G\"{o}del, ``What is Cantor's continuum 
problem?'', {\em
BPPM}, p. 483. The quotation at the beginning of this paragraph is on 
p. 483 of that same
book.} Cantor's Continuum Hypothesis is an important but contested mathematical
statement that brings together countable $\nn$ and uncountable
$2^{\nn}$ infinities in a simultaneous way and with no other cardinality
intervening. It is contested because its validity is undecidable 
within the framework of
Zermelo-Fraenkel axiomatics (equipped with AC) for set theory, as has 
been proven by
G\"{o}del and Cohen. We do not presuppose ZF, however, and will speak 
about sets in an
utterly naive way. We claim that, by generating the 
large[s]-and-small[s], Zeno's
procedure $\mathfrak{Z}$  really does the same as CH. We shall see in 
what follows that
in Zeno's construction do indeed appear the cardinalities associated 
with the naturals
and the reals, and nothing else, so that in the world made up from 
the elements of Zeno's
paradox --- plurality, extension and change ---, Cantor's Continuum 
Hypothesis is valid,
or, put otherwise, that they are aequivalent:

\vspace{-3mm}

\[ \mathfrak{Z} \quad \equiv \quad \aleph_{0} \not= 2^{\aleph_{0}} =
\aleph_{1} \]

{\sc Outline of the proof}.-- It is natural to represent a Zenonian extended
body by a finite continuous line segment.\footnote{VLAS, p. 222, 
shows that of an existant
only extension is considered.} We therefore conceive of
$\mathfrak{Z}$ as applied to such an idealised measuring rod $M$, 
which we  arbitrarily
set equal to unity. Now $M$ will be dichotomised Zenonian-wise, i.e., 
{\bf at infinity,
through and through, and simultaneously}. We label the resulting 
left-hand parts $0$, the
right-hand parts $1$ in every subsequent generation. At each division 
a {\em canonical
choice} is required that sends the divisional sequence down along a 
specific path. We
obtain a collection of nodes and paths, each equipped with a unique binary
representation.\footnote{Weyl was the first to formulate this idea, 
but he linked binary
numbers with the nodes of the tree, thus only capturing the potential 
infinity encoded by
it. Given his finitist a priori this was natural enough. H. Weyl, 
{\em o.c.}, pp. 43-54.
See also Tyles, {\em o.c.}, pp. 64-67.} Thanks to the simultaneity, 
{\em the procedure and
its result coincide,} so that we can look at it constructively, 
without being forced into
an intuitionstic approach. In the
$\nn$-th generation, at the $\omega$-divisional level, the collection 
forms the set $\{0,
1\}^{\nn}$, that contains all possible partitions (nodes) with their parts
(connecting paths). It represents all possible orders on $\nn$. The 
set thus generated
must be aequivalent to the powerset-lattice for $\nn$. But this is 
only a part of the
story, because the pathways followed by generation matter as well. As 
already indicated,
the procedure can be viewed from two different perspectives: a 
stepwise doubling of the
number of cells at every generation, the number of generations being 
$n \leq \omega$, or
the totality of all possible combinations of zeros and ones, viz. the 
power set of ${\cal
P}(\nn) = \{0, 1\}^{\nn}$ {\em with some additional structure}; this 
object being the
simplest possible instance of the {\bf Cantor set}. The first process 
is simultaneously
contained in the second. We thus obtain the two kinds of infinitiy 
constituting a finite
body: an infinity of extended parts --- the {\em ever smaller regions} of the
intuitionists --- which would generate when left on its own an 
infinitely large body. And
an infinity of unextended parts --- points --- which, when left on 
its own, would
generate a body with no extension at all. This absurdity does {\em 
not} arise because
division generates them both together --- {\em kai} --- and at once. 
Conclusion: Zeno's
divisional procedure can be represented by a binary tree, a Zenonian
semi-lattice in which the nodes represent the {\em megala}, while the {\em
micra} are given by its paths, in two different infinities, as follows:\\

\vspace{-2mm}

\SelectTips{cm}{}

\[
\vcenter{\xymatrix{{\Bigl[} \ar@{-}[rrrrrrrr]^-{1 = \bot} &&&&&&&& 
{\Bigr]}\\{\Bigl[}
\ar@{-}[rrrr]^-{0} &&&& {\Bigr]\Bigl[}
\ar@{-}[rrrr]^-{1} &&&& {\Bigr]}\\ {\Bigl[} \ar@{-}[rr]^-{00} && {\Bigr]\Bigl[}
\ar@{-}[rr]^-{01} && {\Bigr]\Bigl[} \ar@{-}[rr]^-{10} && {\Bigr]\Bigl[}
\ar@{-}[rr]^-{11} && {\Bigr]}\\ {\Bigl[} \ar@{-}[r]^-{000} & {\Bigr]\Bigl[}
\ar@{-}[r]^-{001} & {\Bigr]\Bigl[}
\ar@{-}[r]^-{010} & {\Bigr]\Bigl[} \ar@{-}[r]^-{011} & {\Bigr]\Bigl[}
\ar@{-}[r]^-{100} & {\Bigr]\Bigl[} \ar@{-}[r]^-{101} & {\Bigr]\Bigl[}
\ar@{-}[r]^-{110} & {\Bigr]\Bigl[}
\ar@{-}[r]^-{111} & {\Bigr]} }}
\]

\[
\vcenter{\xymatrix@C=.5pc@-1pc{0 & {\quad} & {2^{0}} &{\quad}&& &&&&&&&&&&&&&&
*++[o][F-]{1}
\ar[dllllllll]_-{0} \ar@{.>}[drrrrrrrr]^{\mathbf{1}} \\ 1 & {\quad} & 
{2^{1}} &{\quad}&&
&&&&&&  *++[o][F-]{0}
\ar[dllll]_-{0} \ar[drrrr]^{1} &&&&&&&&&&&&&&&& *++[o][F-]{1}
\ar@{.>}[dllll]_-{\mathbf{0}}
\ar[drrrr]^{1} \\ 2 & {\quad} & {2^{2}} &{\quad}&& && *++[o][F-]{00} 
\ar[dll]_-{0}
\ar[drr]^{1} &&&&&&&& *++[o][F-]{01} \ar[dll]_-{0} \ar[drr]^{1} 
&&&&&&&& *++[o][F-]{10}
\ar[dll]_-{0}
\ar@{.>}[drr]^{\mathbf{1}} &&&&&&&& *++[o][F-]{11} \ar[dll]_-{0} 
\ar[drr]^{1}\\  3 &
{\quad} & {2^{3}} &{\quad}&& *++[o][F-]{000} \ar[dl]_-{0} \ar[dr]^{1} 
&&&& *++[o][F-]{001}
\ar[dl]_-{0} \ar[dr]^{1} &&&& *++[o][F-]{010} \ar[dl]_-{0} \ar[dr]^{1} &&&&
*++[o][F-]{011} \ar[dl]_-{0}
\ar[dr]^{1} &&&& *++[o][F-]{100} \ar[dl]_-{0} \ar[dr]^{1} &&&& *++[o][F-]{101}
\ar@{.>}[dl]_-{\mathbf{0}} \ar[dr]^{1} &&&& *++[o][F-]{110} 
\ar[dl]_-{0} \ar[dr]^{1} &&&&
*++[o][F-]{111} \ar[dl]_-{0} \ar[dr]^{1}\\ {\vdots} & {\quad} & 
{\vdots} &{\quad}&
{\vdots} && {\vdots} && {\vdots} && {\vdots} && {\vdots} && {\vdots} 
&& {\vdots} &&
{\vdots} && {\vdots} && {\vdots} && {\vdots} && {\vdots} && {\vdots} 
&& {\vdots} &&
{\vdots}&& {\vdots}\\ {\omega} & {\quad} & {2^{\nn}} &{\quad} & {\Bigl[}
\ar@{.}[rrrrrrrrrrrrrrrrrrrrrrrrrrrrrr] &&&&&&&&&&&&&&&&&&&& {\bullet}
\ar@{}[r]_{0,1010\dots } &&&&&&&&&&  {\Bigr]}}}
\]
\centerline{The measuring rod $M$ $\{0, 1\}^{\small{\nn}}$, the {\em
Philebian set} \, $\pp$.}\\

\noindent {\em For in two ways it can be said that a distance or a 
period or any other
continuum is infinite [apeiron], viz., with respect to the partitions 
or with respect to
the projecting parts} [Aristotle, {\em Phys.} Z, 2, 263a (24-26)].\\

\vspace{-1mm}

\noindent {\em For whoever divides the continuum into two halves 
thereby confers a double
function on the point of division, for he makes it both a beginning 
and an end} [{\em
Phys.} $\Theta$ \, 8, 263a (23-25)].\\

\newpage

\noindent A formal langiage appropriate to describe such structures 
is {\bf domain
theory}.\footnote{S. Abramski and A. Jung, ``Domain Theory'', in {\em 
Handbook for Logic
in Computer Science}, S. Abramski, D. M. Gabbay and T.S.E. Maibaum 
[eds.], Clarendon
Press, Oxford, 1994, chapters 1,2. I thank Bob Coecke
(Oxford) for clarifying discussions on this subject.} Let us once 
again review some basic
notions. Given is the partially ordered set $(P, \order)$. A subset 
$S$ of $P$ is an {\em
upper set} if $\forall \, s, t \in S$ with
$s \order t : s \in S \Rightarrow t \in S$. The symbol $\uparrow \! s$
\, denotes all elements above an $s\in S$. An element $v$ of $P$ is 
an upper bound for
$S \subseteq P$ if for any $s \in S : \, s \order v$. The least 
element of the set of
upper bounds of $S$ is called the {\em supremum} (or {\em join\,})
$\bigvee \! S$. An element $t \in P$ is {\em maximal} if there is no 
other element of
$P$ above it. A {\em directed upper set or filter}
$S$ of $P$ is a nonempty set of which every pair of elements has an 
upper bound in
$S$. \/If the condition $\uparrow \! s$ \, is met, $S$ is called a 
{\em principal} filter.
The dual notions {\em lower (or down) set} $\downarrow \! s$, lower 
bound, {\em infimum}
(or {\em meet}) $\bigwedge \! S$, minimal element and {\em directed 
lower set or ideal}
appear by reversing the order (arrows). Order in the case of
Zenonian division increases with decreasing intervals, so that $(10)
\pijltje (101)$ implies  $(10) < (101)$. The smallest element is unity, the
measuring rod $M$ itself. Zeno's division operates by the application 
of two simple
rules of mathematical construction, two functions that map partitions on parts:

\begin{enumerate}
\item $r_1 : n \pijlr 2n$ \quad division stepwise \,\, $\leadsto \; 
\text{`small' \,}
{\cal Z}$;

\item $r_2 : n \pijlr 2^n$ \quad simultaneous division, with
${n \order \omega} \,\, \leadsto  \; \text{`large' \,} \mathfrak{Z}$.
\end{enumerate}

\noindent Application of these two rules suffices to build $\mathfrak{Z}$
constructively\footnote{A bit in the spirit intended by
Smyth when he writes: {\em It may be asked whether (...) we adhere to 
constructive
reasoning in our proofs. Actually our procedure is somewhat 
eclectic.} M.B Smyth, {\em
The Constructive Maximal Point Space and Partial Metrizability}, preprint:
http://www.comp.leeds.ac.uk/anthonyr/dtg/papers.htm.}
as the simpliest instance possible of the Cantor set. The 
constructive reasoning used
does {\em not} compell us to an intuitionistic point of view since 
thanks to simultaneity
we can work with actual infinities in a logically consistent way. We 
discussed before the
Zenonian principle of {\bf the aequivalence of the parts and the 
whole.}\footnote{W.E.
Abraham, ``Plurality'', {\em Phronesis}, {\bf 17}, 1972, pp. 40-52.} 
It is essential to
translate this principle, a direct consequence of Zeno's deictically realistic,
exhaustive procedure, into a variant of the basic constructivist 
credo.\footnote{The
importance of an explicit formulation of Zeno's Principle was brought 
home to me during a
discussion with R. Hinnion (ULB), for it clearly demarcates the nature of the
proposed representation. Propositions with regard to CH or AC in the 
following pages are
valid within the confines of Zeno's construction; they do not imply 
any claim with
respect to, say, the Zermelo-Fraenkel axiomatisation of set theory.} 
Let us extract the
principle underlying this construction and summarise it as follows:

\vspace{-2mm}

\begin{flushleft}
\begin{myboxqd}
\centerline{{\bf Zeno's Principle (ZP):} {\em In a construction appear}}
\centerline{{\em only those objects that are constructed effectively,}}
\centerline{{\em by using rules of construction given explicitly.}}
\end{myboxqd}
\end{flushleft}

\noindent The proof basically goes in three steps. Firstly, the 
simultaneous jump from
countably many {\em megala} to uncountably many {\em micra} can be 
rigorously described
as the ideal completion of the small semi-lattice
$({\cal Z}, \order)$ into the large semi-lattice
$(\mathfrak{Z}, \order)$.\footnote{{C.q. cardinal vs. ordinal completion}; W.E.
Abraham, ``Plurality'', pp. 40-52.}  Every node at the
$\omega$-divisional level or, alternatively, every possible branch in the
$\nn^{th}$ generation, i.e., every element $z$ generated by 
$\mathfrak{Z}$, represents a
unique sequence \/ $f(z) = (x_n)_{n \, \in \, \nn}$.\/ The finite ideals
$\downarrow \! x$ can be ordered by inclusion. They all have a 
supremum. The supremum of
the infinite $\omega$-chain of wich they are a part is the maximal 
element $\uparrow \! z
\cap \mathfrak{Z} = \{z\}$. Appropriately enough, a constructive 
interpretation of the
notion of maximality as a criterion or {\em test of fineness} has 
been developed by
Martin-l\"{o}f.\footnote{Discussed in M.B. Smyth, {\em o.c.}, p. 3, 
p. 7 sq.} Each
$\omega$-chain defines a unique order on $\nn$. Every ordinality on 
$\nn$ coincides in a
unique way with such a chain and, by virtue of ZP, no other chains do 
appear. This gives
the orderstructure of ${\cal P}(\nn)$. Therefore the Zenonian 
semi-lattice $(\mathfrak{Z},
\order)$ is a directed and complete partial order or {\em dcpo}.\\

\noindent But $\mathfrak{Z}$ is a total order $\prec$ as well. For by 
the canonical
numbering of the parts, an additional order per generation is 
imposed, and all elements
$z$ will be comparable: $\forall z, z' ; z \prec z' \vee z' \prec z 
\;$;  whence
$(\mathfrak{Z},
\prec)$ is a directed and complete total order or {\em dcto}. This order
catches the influence exercised on the structure by the divisional pathways.
Semantically speaking: the logic encoded by them is intensional, not 
extensional. This
total order is lexicograpic, i.e., according to the {\em principle of first
differences}\footnote{K. Kuratowski and M. Mostowski, {\em Set 
Theory}, North Holland,
Amsterdam, 1968, pp. 224-227.} like in Hausdorff's dictionary.\footnote{R.
Rucker, {\em o.c.}, pp. 82-83.}\hfill$\lhd$\\

\noindent  Our next step is to establish the existence of certain bijective
relationships. This demonstration essentially relies on the
well-known Cantor-Schr\"{o}der-Bernstein theorem. It will not be
admissible, however, to equate the left and right hand variants of the binary
representations of partitions, the points of division at the rational 
multiples of
$2^{-n}$, as in the standard case.\footnote{P.J. Cameron, {\em Sets, Logic an
Categories}, Springer, London etc., 1999, p. 128-129.} They will on 
the contrary be used
as the lefthand $d$ and righthand $d'$ closure of the adjacent, 
non-overlapping parts,
represented by intervals, produced in every generation. This by the 
way completely
justifies Aristotle's seemingly enigmatic comment: {\em For in two 
ways it can be said
that a distance or a period or any other continuum is infinite, viz., 
with respect to the
partitions or with respect to the projecting parts} [{\em Phys.} Z, 
2, 263a (24-26)]. In
honour of Brouwer, we will call them {\bf
\em dubbelpunten} [double points] or simply \, {\bf \em dubbels} \,
[doubles].\footnote{L.E.J. Brouwer, {\em o.c.}, p. 56.}


\[ \; M \, \stackrel{1}\longleftrightarrow \,
\text{semi-lattice} =
\{0, 1\}^{\nn} \, \stackrel{2}\longleftrightarrow \, [0, 1]
\, \stackrel{3}\longleftrightarrow \, \rr \]

\noindent The set of all possible binary
sequences
$f(z)$ is
$\{0, 1\}^{\nn}$. We already saw that $\mid\!\!\{0, 1\}^{\nn}\!\!\mid \, =
2^{\small{\nn}}$. Let us now remove one of every pair of doubles $d$ 
and $d'$ that came
with the rational partitions in
$\mathfrak{Z}$ from our set of sequences, by making a {\bf \em 
canonical choice} for the
lefthand side ($L_c$) or the righthand side ($R_c$) in every instance at every
generation. A canonical choice is a function $f_c \, :  \ A \, \dot{\cup} \, B
\stackrel{\cong}\pijll \, \{0, 1\}^{\nn}$. Let
$d = L_c$. We choose in all cases the zero-side by $L_c$.

\bigskip

\noindent $A = \{(x_n)_n \in \{0, 1\}^{\nn} \ | \ \nexists \, N
\ : \ \forall \, n > N  \ : \  x_n = 1\}$ \\

\noindent $A$ does not contain any sequences that exhibit only ones 
from a certain
$x_n$ on. We filter out the right hand representations of the 
rational multiples of
$2^{-n}$. The so-called redundancy in the binary representation of 
the real numbers is now
removed. Our set $A$ is thus identical with the standard interval
$[0, 1]$. In our case, however, the removed sequences are not simply 
deleted, but
carefully collected in a separate set $B$.

\bigskip

\noindent $B = \{(x_n)_n \in \{0, 1\}^{\nn} \ | \ \exists \, N \ :
\ \forall \, n > N  \ : \  x_n = 1\}$ \\

\noindent This set contains exactly those sequences that do exhibit 
only ones from a
certain $x_n$ onwards; its cardinality will be $\mid\!\!\nn\!\!\mid$. 
 From our definition
it follows that $A = B^c$, such that $A \, \dot{\cup} \, B = \{0, 
1\}^{\nn}$. The
cardinality of the coproduct or disjoint union $A \, \dot{\cup} \, B$ 
is equal to the sum
of the cardinalities of the sets composing it:
$\mid\!\![0, 1]\!\!\mid \, + \, \mid\!\!\nn\!\!\mid = 
2^{\small{\nn}}$. Indeed, from
Cantor's diagonalisation argument we know that $2^{\aleph_0} \, + \,
\aleph_0 = 2^{\aleph_0}$.\footnote{The difficulties going with 
arbitrary sums of infinite
powers are discussed in M. Potter, {\em o.c.}, p. 170 sq.} Equal powers imply
aequivalence; thus we can conclude
$\{0, 1\}^{\nn} \lppijl [0, 1]$. By virtue of its total order, our set $\{0,
1\}^{\small{\nn}}$ is complete, and a complete set which has a 
countably dense subset is
a continuum. And every continuum is isomorphic to the real 
line.\footnote{The conditions
for completeness are either that every non-empty subset of a 
considered set which has an
upper bound has a supremum, or that ever non-empty subset with a 
lower bound has a
minimum. Given the order imposed on $\{0, 1\}^{\small{\nn}}$, the 
second condition is
trivially fulfilled. M. Potter, {\em o.c.}, pp. 119-121.} (The 
aequivalence of $[0, 1]$
with $\rr$ is standard and can be shown by geometrical means.) This 
establishes the
required aequivalences.\hfill$\lhd$\\

\noindent Thirdly regarding order. Our situation reminds us of that 
of infinite sets which, an equal
cardinality notwithstanding, possess different ordinalities. Plato 
makes plain that this
divisional procedure can proceed by ``twoness'', but equally well by 
``threeness'', or
any other number [{\em Philebus}, 16(e)], although the examples given 
are always carried
out by means of {\em aoristos duas}. Therefore we will call in
what follows $\{0, 1\}^{\small{\nn}}$ the {\bf Philebian set}, 
symbolised by $\pp$. This
brings us to the theorem: the sets $\pp$ and
$\rr$ possess equal cardinality but different ordinality\footnote{I 
owe this idea to
a discussion with Tim Van der Linden.}:


\[ \{0, 1\}^{\small{\nn}} \, \text{\em  is \ {\bf \em not 
order-isomorphic} \ to}
\, [0, 1] \]

\noindent Now $X$ is order-isomorphic with $Y$ iff a bijection exists that
satisfies the following condition: $ k \, : \, (X, \order) \pijlr (Y, 
\order)$ with $x
\order y \aequi k(x) \order k(y)$. This implies that, if $(Y,
\order)$ is dense, then $(X, \order)$ will be dense as well. Indeed, 
suppose $Y$ to be
dense. Let $k$ be an order-isomophism, such that $k \ : \ X \pijlr Y 
\ : \ a \, < \, b \,
\aequi k(a) < k(b)$.\\

\noindent Take $c < d \in Y$ then $\exists! \, a, b \in X \, : \, k(a) = c,
\ k(b) = d$.   But $Y$ is dense, therefore $\exists \, m \in Y : \, c 
< m < d$ and so
$a = k^{-1}(c) \leq k^{-1}(m) \leq b = k^{-1}(d)$. However, we 
demonstrated before that
the order generated by $\mathfrak{Z}$ is lexicographic by nature. 
This order is preserved
even for the rational doubles $d$ and $d'$, for which the inequality 
$d < d'$ holds. But
$\nexists \, m \in \{0, 1\}^{\nn} \, : \, d < m < d'$. Then $\{0, 
1\}^{\small{\nn}}$
cannot be dense everywhere. By our canonical choice $L_c$ we removed 
all $d'$ from
$\{0, 1\}^{\small{\nn}}$, thus generating a set $A$ which clearly has 
the property of
being dense, in exactly the same way as $[0, 1]$. Therefore $[0, 1]$ and $\{0,
1\}^{\small{\nn}}$ cannot be order-isomorphic. Furthermore all ideals 
$\downarrow \! x$
included in each unique sequence \/ $f(z) = (x_n)_{n \, \in \,
\nn}$\/ in $\pp$ are finite and non-empty. A well-order $S$ is a
total order in which every non-empty subset of
$S$ has a least element. This is trivially the case here; whence $\pp$ is a
well-order. And well-ordering is aequivalent to the AC. Thus Zeno's
procedure comes about as a precept for the construction of a well-ordered
continuum!\hfill$\lhd$\\

\noindent Now consider the following statement:\\

\noindent {\em Though Zermelo's theorem assures that every set can be 
well-ordered, no
specific construction for well-ordering any uncountable set (say, the 
real numbers) is
known. Furthermore, there are sets for which no specific construction 
of a total order
(let alone a well-order) is known (...)}\footnote{J. Dugundji, {\em 
Topology}, Allyn and
Bacon, Boston, 1966, p. 35.}\\

\noindent {\sc Conclusion}.-- For Zeno's divisional procedure we proved the
following\footnote{C. Vidal (Sorbonne) attracted my attention to a 
book by A.W. Moore,
{\em The Infinite}, Routledge, London and N.Y., 2001 [1990], in 
which, though without
invoking Zeno, a promising line of ideas is developed, proposing a link between
``infinity by addition'' and ``infinity by division'' --- clearly the 
same concepts as
our ``stepwise'' and ``simultaneous'' division ---, the continuum 
hypothesis and the
cardinalities of $\nn$ and $\rr$, but only to dismiss the 
possibility! The reason
apparently is that it implies according to Moore a variant of the 
Bernoullian fallacy,
given his reference to (rational) density. See pp. 154-158 of the 
cited work. I want to
stress once more that Zeno does not commit any such fallacy, and 
refer to our own, but
also to Abraham's analysis of Zeno's ideas. For the latter, see 
Abraham's already
referenced Plurality-article.}\\

$\star$ \, Divisional procedure \, $\mathfrak{Z} \equiv 
(\mathfrak{Z}, \prec)$ \quad
\quad (simultaneity)

$\star$ \, Immediate successor \quad (ideal completion)

$\star$ \, Two different kinds of infinity	$\mid\!\!\nn\!\!\mid 
\, \not= \,
\mid\!\!2^{\footnotesize{\nn}}\!\!\mid$ (Cantor's Theorem)

$\star$ \, Axiom of Choice \quad (well-ordering)

$\star \,\, M \, \stackrel{1}\longleftrightarrow \,
\text{lattice} =
\{0, 1\}^{\nn} \, \stackrel{2}\longleftrightarrow \, [0, 1]
\, \stackrel{3}\longleftrightarrow \, \rr$	\quad (Cantor-Bernstein)\\

\vspace{-4mm}

\begin{flushleft}
\begin{myboxqd}
\vspace{1mm}
\noindent Contrary to the above statement, our claim is that Zeno's 
divisional procedure
provides a specific way to construct a well-ordered continuum 
isomorphic to the real
number line (though {\em not} \/to $\rr$), in which Cantor's 
continuum hypothesis is
valid:

\vspace{-5mm}

\[ \mathfrak{Z} \quad \equiv \quad \aleph_{0} \not=
2^{\aleph_{0}} =
\aleph_{1} \]\vspace{-3mm}
\end{myboxqd}
\end{flushleft}

\vspace{5mm}

{\sc Remark}.-- Every real number in this construction is represented 
and defined
by a unique order on the set of natural numbers, a specific subset $f(z) \in
{\cal P}(\nn)$. Moreover, thanks to ZP all orders on $\nn$ are 
included, and nothing else,
thence the argument is truly constructive. The ``line'' formed by the 
divisional {\em
loci} $(d_i,d_i')$ is isomorphic to
the rational line. Is it possible to get the rational line per se, so 
that the point
intervals to which each of the $(d_i,d_i')$ supplies an ending and an 
opening bracket
contain exactly the irrationals? A promising line of thought could be 
to generalise the
divisional procedure over all prime numbers, a suggestion bestowed 
with some authority by
tradition.\footnote{J. Stenzel, {\em o.c.}, p. 53 sq. refers to a 
passage in Book $A$ of
the {\em Metaphysics}. See especially p. 56.} It is clear that one 
should include as well
the powers of the primes $n$p$^k$, with $n, k \in \nn$. {\em All}\/
rationals are doubly represented in that case, an apt way to enhance their
character of finite fractions. But even then not all possible 
intervals are covered,
thanks to the one domensional a priori implicit in the division. This 
not yet suffices to
recover all ordinals of ordertype zero, for we still do not cover all possible
divisional intervals. In order to do that, one should introduces 
multiples of powers of
primes. One then of course calls upon the availability of the natural 
numers, but this is
no problem, since $\mathfrak{Z}$ generates these. This would, 
however, constitute a
deviation of Zeno's intentions --- the one-dimensional a priori being 
inhaerent to his
procedure --- which is neither appropriate, neither necessary here. 
So we are left with an
infinity of infitely tiny holes in Zeno's strange continuum. 
Intuitively, they are so
small as to be numbers smaller than whatever real number given, but 
different from
$0: h < r \ \text{with} \ r \in \rr$. In other words, they do not 
obey the Archemidean
axiom: $\forall \, x, y \in \rr$ with $x < y \, : \, \exists \, m \in 
\nn \, : \, m.x >
y$, which precisely excludes the existence of infinitesimally small 
numbers.\footnote{ R.
Courant and F. John, {\em Introduction to Calculus and Analysis}, vol. I,
Wiley/Interscience, 1965,, p. 94.} Put otherwise: $h$ is an 
infinitesimal iff \/$\forall
\, m \in \nn \ : \ \M m.h \M \, < 1\,$. This implies that, whatever 
they are, they are
not reals. But there are numbers in mathematics that fulfill the 
criterium of being
smaller than whatever given real. They are called {\bf 
infinitesimals}. They suffer from
a bad reputation, however, for they are held responsible for the 
notoriously shaky nature
of the foundations on which the early calculus rests. It is well 
known that the {\em locus
classicus} of modern natural science, Newton's monumental {\em Principia}
[1687]\footnote{A modern edition known as the {\em variorum} edition 
has been edited by
A. Koyr\'{e} and I.B. Cohen, {\em Isaac Newton's Philosophiae 
Naturalis Principia
Mathematica, The Third Edition (1726) with variant readings}, 
Cambridge University Press,
Cambridge, 1972.}, uses a theory of infinitesimals in a geometrically 
disguised manner.
The underlying {\em theory of fluxions}, which explicitly uses 
infinitesimals, got
published only afterwards, although it was developed twenty years 
before the {\em
Principia} appeared. A clear exposition of the fluxion-theory cast in 
a geometrical
framework can be found in the short tract {\em De Quadratura 
curvarum}, published as an
appendix to Clarke's latin translation of the {\em Opticks} 
[1706].\footnote{A discussion
of Newton's underlying `finitist' approach towards descriptions of 
natural phaenomena, as
well as his reasons for sticking to the classic geometrical approach, can be
found in G. Guicciardini, {\em Reading the Principia. The Debate on 
Newton's Mathematical
Methods for Natural Philosophy from 1687 to 1736}, Cambridge 
University Press, Cambridge,
1999/2003, p. 27 sq.} That other giant, Leibniz, published his own 
version of the calculus
as the {\em Nova Methodus} in the {\em Acta Eruditorum} 
[1684].\footnote{G.W. Leibniz,
``Nova methodus pro maximis et minimis, itemque tangentibus, quae nec 
fractas nec
irrationales quantitates moratur, et singulare pro illis calculi 
genus'', {\em Acta
Eruditorum}, vol. III, 1684.} The notations currently in use in 
calculus are introduced
therein, together with a proof of Leibniz's ``chain rule''. We 
already explained that an
infinitesimal number is a number smaller than any given real number, 
while it remains
different from $0$. Newton and Leibniz dealt with these quantities in 
their attempts to
formalise consistently the apparently natural notions shoring up the 
newly emerging
infinitesimal calculus. But they hit upon deeply rooted logical paradoxes.
Berkeley with apparent taste exposed them in a devastating way in his {\em The
analyst}.\footnote{I consulted the electronic edition edited by D.R. Wilkins:
http://www.maths.tcd.ie/pub/HistMath/People/\\Berkeley/Analyst/Analyst.html.} 
The problem
essentially boils down to the fact that something being in the sense 
of extended)
transforms into something being not, whence his notorious expression 
that they are {\em
the ghosts of departed quantities}. These logical problems 
concomittant to the direct use
of infinitesimals in calculus led to the development of the Cauchy 
theory of limits, and
the subsequent reformulation of the basic definitions of calculus in terms of
Weierstrassian $\epsilon$--$\delta$ formalism.\footnote{R. Courant 
and F. John, {\em
o.c.}, vol. I, pp. 95-97.} One can however doubt whether the problem 
really dissappeared
because an axiom of continuity has to be invoked to render the method 
unambiguous. The
concept of infinitesimal resurfaced again during the second half of 
the twentieth
century, thanks to the work of Abraham Robinson, who explicitly 
considers his work as an
actualisation of Leibniz's ideas. In order to avoid the logical 
problems intrinsic to
them, he introduced infinitesimals defined as {\em hyperreal numbers} 
by means of {\em
model theory}.\footnote{A. Robinson, {\em Non-Standard Analysis}, 
North-Holland,
Amsterdam/London, 1966/1974. For Leibniz, see p. 2; p. 260 sq. Also 
J.M. Henle and E.M.
Kleinberg, {\em Infinitesimal Calculus}, Dover, 1979/2003.} In this paper, a
minimalistic, at first purely algebraic approach will be followed 
based on {\em dual}
or {\em nilpotent} numbers, in which the {\em d\'{e}tour} via model 
theory can be
avoided. It is an alternative to a more geometrical approach developed by
Bell.\footnote{Nilpotent infinitesimals have been re-introduced by 
J.L. Bell in the realm
of intuitionistic mathematics. The standard reference is: J.L. Bell, 
{\em A Primer of
Infinitesimal Analysis}, Cambridge University Press, Cambridge, 
1998.} It will bring
about its own logical problems nevertheless, while being entitled to 
a historical
provenance of an ancestry comparable to the other systems. Now the 
suggestion is to fill
out the divisional holes by means of infinitesimals. Let us follow up 
this suggestion.\\

{\sc Zenonian Infinitesimals}.-- Although the set $\pp$ generated by
$\mathfrak{Z}$ clearly is a continuum, we would rather like to be 
certain with regard to
its precise nature, for we would like to do analysis and other things 
that allow us to
apply mathematics on the real world, as in the case of the paradoxes 
of motion. On the
other hand, our problem might at first not seem very serious, for it 
apparently suffices
to systematically equate the unequal objects $d_i$ and
$d'_i$ to get back to $\rr$, and be able to do calculus perfectly 
well. But that would
bring us off the track we set out at the start, which is to stay 
loyal to Zeno's
intentions as much as possible. Luckily enough a way out has been 
opened up for us by
Henri Poincar\'{e}. It consists of carefully distinguishing between 
the different
notions of ``continuum'' which in general are uncritically mixed up, 
and to make explicit
these differences formally. Developing an idea he already launched in {\em La
science et l'hypoth\`{e}se}\footnote{{\em SHYP},
p. 51 sq.}, Poincar\'{e} writes: {\em Il arrive que nous sommes 
capables de distinguer
deux impressions l'une de l'autre, tandis que nous ne saurions 
distinguer chacune d'elles
d'une m\^{e}me troisi\`{e}me.}\footnote{SVAL, p. 61.} We could, say, 
distinguish a weight
of 12 grammes from one of 10 grammes, while the intermediate weight 
of 11 grammes would be
indistinguishable from either of both.  In our experience of physical 
reality, there
would be a continuum between them. This amounts into the following paradoxical 
definition
of {\em le continu physique}:

\[ \text{A = B, B = C, A $<$ C} \]

\noindent which clearly violates the {\em Principium 
Contradictionis}. But suppose we
enlarge our perceptive capacities; would then the difficulty not 
simply disappear? No,
for it would be easy enough to find elements D, E that can be 
intercalated between A, B
so that

\[ \text{A = D, D = B, A $<$ B; B = E, E = C, B $<$ C} \]

\noindent and so on, ad infinitum, in analogy with infinite division. The
difficulty does only recede, not disappear. It is remarkable how close this
comes to Brouwer's original starting point: {\em we [zullen] nader ingaan op de
oer-intu\"{\i}tie der wiskunde (...) als het van qualiteit ontdane 
substraat van alle
waarneming van verandering, een eenheid van continu en discreet, een 
mogelijkheid van
samendenken van meerdere eenheden, verbonden door een ``tussen'', dat 
door inschakeling
van nieuwe eenheden, zich nooit uitput.}\footnote{{\em we will take a 
closer look at the
fundamental mathematical intuition (...) as the substrate to all 
perceptions of change,
stripped off all its qualities, a unity of continuous and discrete, a 
possibility of the
thinking-together of several individualities linked by a ``between'' 
that will never be
exhausted by adding new individualities in between.} L.E.J. Brouwer, 
{\em o.c.}, p. 43.
[My translation.]} This is why Poincar\'{e} calls the physical 
continuum {\em une
n\'{e}buleuse non r\'{e}solue.} The mathematical continuum will serve 
to resolve this
cloud and thus to remove the intolerable contradiction: {\em c'est le continue
math\'{e}matique qui est la n\'{e}buleuse r\'{e}solue en 
\'{e}toiles.}\footnote{{\em
SVAL}, p. 61.}  That this description really catches the notion of 
division we explicated
in $\mathfrak{Z}$ on behalf of Zeno is plain, because of the decisive 
property it shares
with the latter: {\em Celui-ci est une
\'{e}chelle dont les \'{e}chelons (nombres commensurables ou
incommensurables\footnote{i.e., rational or irrational numbers!}) 
sont en nombre infini,
mais {\bf \em sont ext\'{e}rieurs les un aux autres, au lieu 
d'empi\'{e}ter les un sur les
autres} comme le font, conform\'{e}ment
\`{a} la formule pr\'{e}c\'{e}dente, les \'{e}l\'{e}ments du continu 
physique.} Now let
us give the name of {\bf Poincar\'{e} continuity} to this paradoxical 
property of the
physical continuum. I remind the reader here of Aristotle's utterly 
correct observation:
{\em For whoever divides the continuum into two halves thereby 
confers a double function
on the point of division, for he makes it both a beginning and an 
end.} Our claim with
respect to the  difference between $\rr$ and $\pp$ can be summarised as:

\[ \text{\em {the set $\rr$ is Poincar\'{e} continuous,}} \]

\noindent and therefore represents the {\em physical} instead of the 
mathematical
continuum, contrary to the standard view. It is the set  $\{0, 
1\}^{\small{\nn}}$ which,
by appropriately discriminating between $d$ and $d'$, resolves the 
unresolved nebula and
thus represents the {\em mathematical} continuum, precisely because 
it contains an
infinity of infinitely tiny holes. We showed however that this is
indeed the case, and will now try to find a way to render the 
mathematical continuum
intuitively continuous again.\\

We look for a notion of infinitesimals that does not treat them as the mere
logical consequence of the introduction of a special operator, as is 
the case with
Robinson (although we retain his notion of `hyperreality'), but as 
mathematical entities
in their own right. Infinitesimals of this kind have been introduced 
at the very outset of
the development of modern natural science and the calculus by {\bf 
Bernard Nieuwentijt}, a
seventeenth century Dutch mathematician and theologian.\footnote{I 
encountered this
name for the first time in an article by J. L. Bell, ``Infinitesimals and the
Continuum'', {\em Mathematical Intelligencer} , {\bf 17}, 2, 1995, 
ftn 2. A more
elaborate treatment in J.L. Bell, {\em The Continuous and the Infinitesimal in
Mathematics and Philosophy}, Polimetrica, Milano, 2005.} He defends 
them in a controversy
with Leibniz over the years 1694-1696 about presumed inconsistencies 
in the latter's
variant of the calculus. Nieuwentijt's main criticism obviously 
concerns the use of
higher order differentials.\footnote{R.H. Vermij, {\em Secularisering 
en natuurwetenschap
in de zeventiende  en achtiende eeuw: Bernard Nieuwetijt}, Rodopi, 
Amsterdam, 1991, p. 24
sq.; P. Mancosu, {\em Philosophy of mathematics and mathematical 
practice in the
seventeeth century}, Oxford University Press, N.Y., 1996.} It thus 
comes close to
Berkeley's attack on the Newtonians in his {\em Analyst}\footnote{N. 
Guicciardini, {\em
o.c.}, pp. 199-200.}, though on radically different grounds, for he 
does not reject ---
as with Berkeley --- Newton's use of the potential infinite; he on 
the contrary accepts
(on logical and theological grounds) the actual infinite as well! 
Nieuwentijt's approach
has --- in certain respects --- a very `modern' smell about it. He 
was unsatified by the
{\em ad hoc} nature of many solutions to specific problems presented 
by Newton and
Leibniz. He pursued to develop a logically consistent, general 
foundation for the
calculus, based on a study of the properties of the infinite. This 
problem is related to
the ontological status of infinitesimal quantities: do they really 
exist or are they
merely limiting cases, i.e., finite approximations? It will be clear 
that in the latter
case only potential infinity is assumed, in accordance with 
Archimedes's geometrical
method of exhaustion, to which both Newton and Leibniz painstakingly 
show their methods to
correspond.\footnote{N. Guicciardini, {\em o.c.}, p. 164 sq.} 
Nieuwentijt's criticism is
that such an approximative description cannot supplant an exact 
definition, hence cannot
serve as a basis for valid logical deductions in calculus. He 
formulates a fundamental
axiom on the basis of which infinitesimals could be exactly defined: 
{\em everything
which, when multiplicated by an infinitely great number, does not 
render a given [finite]
number, whatever small, cannot be counted as a being but should in 
the realm of geometry
be considered as a pure nothing}. On this basis he formulates the 
arithmetical rules for
the infinite. He moreover proves that from this axiom it follows that 
the square of an
infinitesimal --- taken as the number $a = A/m$, with $m = \infty$ 
and $A$ finite  ---
must be zero. Indeed, multiplying $a^2$ by $m$ is by definition 
$(\frac{A}{m})^2m$, and
this equals $\frac{A^2}{m}$, again an infinitesimal, and so by the 
grounding axiom equal
to $0$.\footnote{R.H. Vermij, {\em o.c.}, p. 19. A more detailed exposition of
Nieuwentijt's mathematical methods in H. Weissenborn, {\em Die 
Prinzipien der h\"{o}heren
Analysis in ihrer Entwicklung von Leibniz bis auf Lagrange, als ein 
historisch-kritischer
Beitrag zur Geschichte der Mathematik}, H.W. Schmidt, Halle, 1856.} 
For whatever you do
to get it over the border of mathematical visibility will fail, and 
so it must be
no-thing.  Nieuwentijt thence defines infinitesimals as follows: {\em 
Si pars qualibet
data minor $b/m$ ducatur in se ipsam, vel aliam qualibet data minorem
$c/m$, erit productum $bb/mm$ seu $bc/mm$ aequale nihilo seu non 
quantum. [When a certain
infinitesimal part $b/m$ is applied onto itself, or on another 
infinitesimal part $c/m$,
the product $bb/mm$ or $bc/mm$ will be equal to zero, or have no
quantity.]}\footnote{Bernhardi Nieuwentijt, {\em Analysis infinitorum 
seu curvilineorum
proprietates ex polygonorum natura decuctae}, Amstelaedami, Wolters, 
1695, Praefatio,
lemma 10. That he really has infinitesimals in mind becomes clear 
from the definitions on
p. 1 of his book: ``data major'' = ``infinitam'' and ``data minor'' = 
``infinitesimam''.
Thus te bigger $m$ becomes, the smaller $b/m$ will be.} Although 
themselves different from
zero, the infinitesimals are so small that their square vanishes. 
Intuitively this is not
so strange as it may seem at first glance: consider the square of the 
real number
$0, 0001$! Thus the infinitesimal character becomes a definable 
property in its own
right. This is where the concept ``nilpotency'' comes in.\\

Now let us look in more detail at the set of Nieuwentijt
infinitesimals $\bb$ (`B' from Bernard) of nilpotent numbers.\footnote{The
initial suggestion that dual or nilpotent infinitesimals are the best 
match for our
purposes stems from Didier Deses (VUB). Concerning the algebra, W. 
Lowen (Paris VII) made
many valuable suggestions.} We consider
$\bb$ as a plane
$\rr
\times
\rr$ in which the number $(a, b)$ consists of a real component $a$ 
and a hyperreal
component
$b$. The latter one is the coefficient of the infinitesimal $h$. 
Addition of two such
hyperreals is component-wise: $(a + bh) + (c + dh) = (a + b) + (c + d)h$. For
multiplication we have: $(a + bh) \cdot (c + dh) = ac + adh + bhc + bdh^2
= ac + (ad + bc)\;h$, or $(a,b)(c,d) = (ac, ad + bc)$. In other 
words, we multiplicate $a
+ bh$ en $c + dh$ \/as polynomes in the variable $h$, with $h^2 = 0$. 
This of course holds
true for all multiples of $h^2$. As we noted already, our number $(0, 
1)$ is so small
that,  although itself different from $0$, its square equals 
$0$.\footnote{Division can be
defined by means of the adjoint. It is the inverse operation of 
multiplication, with the
{\em caveat} that, in order to remain consistent, $b \not= 0$.}

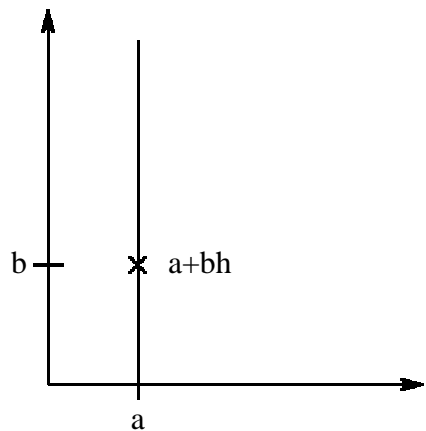
\begin{floatingfigure}{6.5cm}
\setlength{\unitlength}{0.00087489in}
\begingroup\makeatletter\ifx\SetFigFont\undefined%
\gdef\SetFigFont#1#2#3#4#5{%
       \reset@font\fontsize{#1}{#2pt}%
       \fontfamily{#3}\fontseries{#4}\fontshape{#5}%
       \selectfont}%
\fi\endgroup
\begin{picture}(2512,2572)(0,-10)
\thicklines
\drawline(240,285)(2490,285)
\drawline(2370.000,255.000)(2490.000,285.000)(2370.000,315.000)(2406.000,285.000)(2370.000,255.000)
\thicklines
\drawline(240,285)(240,2535)
\drawline(270.000,2415.000)(240.000,2535.000)(210.000,2415.000)(240.000,2451.000)(270.000,2415.000)
\drawline(780,2355)(780,285)
\thicklines
\drawline(150,1005)(330,1005)
\drawline(780,375)(780,195)
\drawline(825,1050)(735,960)
\drawline(735,1050)(825,960)
\put(735,15){\makebox(0,0)[lb]{\smash{{\SetFigFont{12}{14.4}{\rmdefault}{\mddefault}{\updefault}a}}}}
\put(15,960){\makebox(0,0)[lb]{\smash{{\SetFigFont{12}{14.4}{\rmdefault}{\mddefault}{\updefault}b}}}}
\put(960,960){\makebox(0,0)[lb]{\smash{{\SetFigFont{12}{14.4}{\rmdefault}{\mddefault}{\updefault}a+bh}}}}
\end{picture}
\caption{\label{Figuur-Duaal-Getal}The nilpotent number $a+bh$}
\vspace{2mm}
\end{floatingfigure}

\noindent Thanks to these
properties, $\bb$ is a commutative ring with unity, which forms an 
algebra over the
field $\rr$. This $\rr$-algebra is the ringtheoretic quotient $\rr[h]/(h^2)$
of $\rr[h]$. $\rr[h]$ is the polynome ring in the variable $h$,
where $(h^2)$, the ideal generated by $h^2$, is divided out, so that the ring
$\rr[h]/(h^2)$ really is the ring of nilpotent numbers. The 
representation in $\rr
\times \rr$ allows for a geometrical interpretation of $h_r = a + bh$ for
$r \in \rr$, as can be seen in figure 1. This figure gives further 
information on the
ordertopology on $\bb$. The hyperreal part $b$ clearly assigns a 
unique position to $bh$
on the straight line parallel to the ordinal through $a$. These 
parallels in their turn
occupy a unique position on the abscissa, assigned by $a$. The 
relation satisfies the
criterion for total order: $\forall \, h_r \in \bb : (a, b) \order 
(c, d)$ or $(c, d)
\order (a, b)$. The total order on $\bb$ will be determined by means 
of the {\em first
difference}, so that it is not merely total, but lexicographical as 
well. This is in
agreement with what we obtained when dividing Zenonian-wise. It 
furthermore will be
possible to embed $\rr$ in $\bb$ in by means of the injection $\iota 
: \rr \krulpijlr
\rr \times \rr$ : $r \mapsto h_r$ with $h_r = (r, 0)$. This embedding 
preserves the order.
Let us now return to our initial question: given that the 
lexicographical order that
equally governs $\{0, 1\}^{\small{\nn}}$ prevents the double points 
$d < d'$ from
coinciding, how to fill out the gaps in $\pp$ that arise as a 
consequence of their
presence? Our answer will obviously be that the Nieuwentijt 
infinitesimals lie exactly
between them. But this answer implies that it should be possible to 
somehow  construct a
viable completion of
$\pp \times \pp$ in $\bb$, and show that the result remains aequivalent to
$\mathfrak{Z}$. We are facing immediately a problem here, because, 
while the number of
points $a$ on the abscissa is uncountable, we only have a countable 
number of places
$(d,d')$ available to insert what I propose to call a {\em prime 
needle} for reasons that
do not concern us here. But let us neglect this problem at first by arbitrarily
supposing that not every $a$ possesses this power to instantiate a hyperreal
monad\footnote{This terminology stems from Robinson's seminal work, 
and, although a bit
unfashionable these days, I insist on using it as a tribute to him. 
Cfr. A. Robinson,
{\em o.c.}, p. 57.}, and look what happens when we insert a prime needle
$D_i$ between the members of each couple
$(d_i, d'_i)$, thus executing graphically the construction exactly as 
we need it.
This construction does indeed generate the Euclidean plane, thanks to an
at first glance improbable {\sc theorem} S proven by W. 
Sierpi\'{n}ski that states: \[
\text{\em {the plane is a sum of a countably infinite number of curves.}}\]

\begin{figure}[h]
\vspace{1mm}
\setlength{\unitlength}{0.00087489in}
\begingroup\makeatletter\ifx\SetFigFont\undefined%
\gdef\SetFigFont#1#2#3#4#5{%
      \reset@font\fontsize{#1}{#2pt}%
      \fontfamily{#3}\fontseries{#4}\fontshape{#5}%
      \selectfont}%
\fi\endgroup
\begin{picture}(5254,2389)(0,-10)
\drawline(552,192)(552,12)
\drawline(462,102)(102,102)
\drawline(12,192)(12,12)
\drawline(822,192)(732,192)
\drawline(732,192)(732,12)
\drawline(732,12)(822,12)
\drawline(822,102)(1182,102)
\drawline(1182,192)(1272,192)
\drawline(1272,192)(1272,12)
\drawline(1272,12)(1182,12)
\drawline(102,192)(12,192)
\drawline(12,12)(102,12)
\drawline(462,192)(552,192)
\drawline(552,12)(462,12)
\drawline(642,102)(642,2172)
\put(462,282){\makebox(0,0)[lb]{\smash{{\SetFigFont{12}{14.4}{\rmdefault}{\mddefault}{\updefault}d}}}}
\put(732,282){\makebox(0,0)[lb]{\smash{{\SetFigFont{12}{14.4}{\rmdefault}{\mddefault}{\updefault}d'}}}}
\thicklines
\drawline(2982,102)(5232,102)
\drawline(5112.000,72.000)(5232.000,102.000)(5112.000,132.000)(5148.000,102.000)(5112.000,72.000)
\thicklines
\drawline(2982,102)(2982,2352)
\drawline(3012.000,2232.000)(2982.000,2352.000)(2952.000,2232.000)(2982.000,2268.000)(3012.000,2232.000)
\drawline(3162,2172)(3162,102)
\drawline(3342,2172)(3342,102)
\drawline(3522,2172)(3522,102)
\drawline(3702,2172)(3702,102)
\drawline(3882,2172)(3882,102)
\drawline(4062,2172)(4062,102)
\drawline(4242,2172)(4242,102)
\drawline(4422,2172)(4422,102)
\drawline(4602,2172)(4602,102)
\drawline(4782,2172)(4782,102)
\drawline(4962,2172)(4962,102)
\drawline(2982,282)(5052,282)
\drawline(5052,462)(2982,462)
\drawline(5052,642)(2982,642)
\drawline(5052,822)(2982,822)
\drawline(5052,1002)(2982,1002)
\drawline(5052,1182)(2982,1182)
\drawline(5052,1362)(2982,1362)
\drawline(5052,1542)(2982,1542)
\drawline(5052,1722)(2982,1722)
\drawline(5052,1902)(2982,1902)
\drawline(5052,2082)(2982,2082)
\end{picture}
\caption{\label{Figuur-HR-Getal}The nilpotent numbers complete $\pp 
\times \pp$ into
$\bb$}
\end{figure}
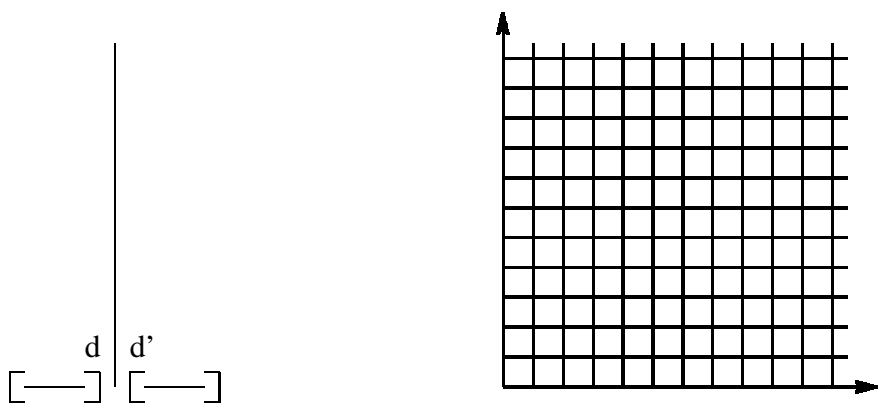

\noindent Since we demonstrated before the aequivalence of 
$\mathfrak{Z}$ with CH, for our
needs aequivalence of S with the latter would do. But this is exactly 
Sierpi\'{n}ski's
point: to prove that theorem S quoted above is aequivalent to CH.\footnote{W.
Sierpi\'{n}ski, {\em Hypoth\`{e}se du Continu}, Z Subwencji Funduszu 
Kultury Narodowej,
Warsawa/Lw\'{o}w, 1934, pp. 9-12. [Reprinted by Chelsea Publishing 
Company 1956]}
Sierpi\'{n}ski's starts by proving the weaker {\sc theorem} 
S$^{\ast}$: {\em The set of
all points of the plane is itself a sum of two sets of which one is 
at most countable on
every straight line parallel to the ordinate, the other one at most 
countable on every
parallel to the abscissa.} Our initial hypothesis that only on the 
rational number of
{\em loci}
$(d_i,d_i')$ the real part $a$ will possess the power to intantiate a 
hyperreal monad
$bh$ along the prime needle $D_i$,  Our embedding now of course will 
be $\epsilon : \pp
\krulpijlr \pp \times \pp$ : $r \mapsto h_r$ with $h_r = (r, 0)$, 
which preserves
the order. But it is far from clear yet what this
exactly means in terms of the {\em Gedanken}-experiment which is at 
the heart of Zeno's
procedure. Therefore we will have to elucidate further its physical
meaning, by taking a closer look at its geometrical implications. At 
first we need to
discuss the paradoxes of motion, and see how tehy fit into the play.\\

{\sc Zenonian plurality in the Motion Paradoxes}.-- The classical
interpretation again presents us Zeno's argumentative strategy as a 
dilemma, at bottom
concerning the continuity of space, and, as a corrolary, the continuity of
time.\footnote{W. C., Salmon (ed.), {\em Zeno's Paradoxes}, Hackett, 
Indianapolis, 2001
[1970]. Owen succeeds in bringing all motion paradoxes on a par by referring
explicitly to the plurality problem. G.E.L. Owen, {\em op. cit.}} An
alterntive though related interpretation centers upon the denseness 
postulate for both
space and time.\footnote{A. Gr\"{u}nbaum, {\em o.c.}, p. 37 sq.} We 
will shortly consider
the applicability of our representation
$\mathfrak{Z}$ to Zeno's motion paradoxes, confining ourselves to an 
outline of the
possibility, and leaving a more elaborate discussion for future work. 
It will be
indicated furthermore how this links to the question at the end of 
the preceding
paragraph. This comes down to devising an interpretation that subsumes PM under
Zeno's simultaneous and through and through divisional procedure (as 
we modelled it),
which indeed matches the nature of our claim with respect to them.

Zeno's famous Paradoxes of Motion are transmitted to us by Aristotle [{\em
Phys.}, Z 9, 239b], with the comment that they are notoriously difficult to
refute. And indeed attempts to either refute, either resolve them 
have been at the order
of the day up to the present: {\em no one has ever touched
Zeno without refuting him, and every century thinks it worthwhile to refute
him}.\footnote{A.N. Whitehead, {\em Essays in Science and 
Philosophy}, Philosophical
Library, New York, 1947.} Let it suffice to say that, however 
relevant in themselves for
future developments in, say, mathematics, {\em all} presumed 
refutations hinge on
non-Zenonian praemisses, so that, whatever it was that was refuted, 
it was certainly not
Zeno. Again, it is not Zeno who presupposes space or time; nor does he assume
hypothetical postions on their being discrete or continuous. The only 
thing needed to
find PM back is Zeno's simultaneous and through and through division. 
Our analysis of PP
thus should be qpplicable to PM as well.  The reason invariably is 
that ``to move"
implies ``to count the uncountable'', or, which boils down to the 
same, to measure
implies to apply commensurable units to incommensurable quantities.

\SelectTips{cm}{}

\begin{figure}[h]
\vspace{-1mm}
\[
\vcenter{\xymatrix{ {\Bigl[} \ar@{-}[r]^-{000} & {\Bigr]\Bigl[}
\ar@{-}[r]^-{001} &  {\Bigr]\Bigl[} \ar@{-}[r]^-{010} & {\Bigr]\Bigl[}
\ar@{-}[r]^-{011} &  {\Bigr]\Bigl[} \ar@{-}[r]^-{100} & {\Bigr]\Bigl[}
\ar@{-}[r]^-{101} &  {\Bigr]\Bigl[} \ar@{-}[r]^-{110} & {\Bigr]\Bigl[}
\ar@{-}[r]^-{111} &  {\Bigr]}\\ {\Bigl[} \ar@{-}[rrrr]^-{0} &&&& {\Bigr]\Bigl[}
\ar@{-}[rr]^-{10} && {\Bigr]\Bigl[} \ar@{-}[r]^-{110} &  {\Bigr]\Bigl[}
\ar@{-}[r]^-{111} & {\Bigr]}\\ {\Bigl[} \ar@{-}[r]^-{000} & {\Bigr]\Bigl[}
\ar@{-}[r]^-{001} &  {\Bigr]\Bigl[} \ar@{-}[rr]^-{01} && {\Bigr]\Bigl[}
\ar@{-}[rrrr]^-{1}  &&&& {\Bigr]} }}
\]
\caption{\label{Figuur-Bewegingsparadoxen}The paradoxes of motion}
\end{figure}
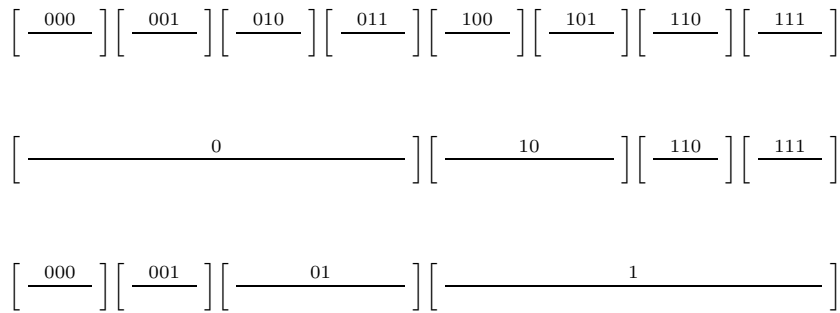

\noindent Bell puts it like this: {\em Continuity and discreteness 
are united in the
process of measurement in which the continuous is expressed in terms 
of separate units,
that is numbers.}\footnote{J.L. Bell, {\em Oppositions and Paradoxes 
in Mathematics and
Philosophy}, [forthcoming in Axiomathes].} This is the fundamental 
point one looses out
of sight when one throws actual infinity out. One assumption 
implicitly underlying the
execution of Zeno's procedure --- we touched it already before, but 
here it becomes
particularly relevant --- is the one-dimensional, horizontal orientation of the
divisional process. This assumption allows Zeno to go over without 
any further ado from
plurality at rest to plurality in motion. When one reads the 
PM-fragments from this
angle, the astonishing coherence of the arguments immediately hits 
the eye. The problems
treated by Zeno as embodied in the four arguments given by Aristotle 
are, by levels of
ascending perplexity\footnote{Traditionally the numbering of the two 
last ones is
switched, although it is clear that the Arrow concludes Zeno's 
stupendous argumentative
sequence. We concord partially with Owen's scheme.}:

\begin{tabbing}
	\= i)  motion cannot take a start;\\
	\> ii)  motion, once started, cannot be completed;  \\
	\> iii) moving bodies passing each other are subject to the 
same paradoxes; \\
	\> iv)  motion is self-contradictory.
\end{tabbing}

\noindent Zeno's point precisely is to show that, {\em however small 
the distance}, the
number of parts to cross will remain the same, i.e. $2^{\nn}$, while 
the number of
partitions (steps) will never exceed $\nn$, and in reality be only 
finite. Time, being
related to distance, is irrelevant in exactly the same way as the 
length of our measuring
rod was irrelevant. It would moreover be quite strange that Zeno, in 
order to defend
Parmenides's stance with respect to the deictical unreality of time, 
would introduce it
to make his point. `To count the uncountable' is thus the motion-face of the
plurality-coin, which, as the reader will remember, can be summed up 
in the slogan `to
consist out of parts with and without magnitude'. One sees that 
Aristotle's choice, far
from being arbitrary, was to pick out exactly those renderings (from 
an undoubtedly
larger corpus) that develop the paradox step by step, in order to lay 
bare its many
faces, and to bring out why it is so difficult to resolve. But 
although the Stagirite
realised the nature of the underlying problem, he apparently did not 
believe that Zeno
himself did. This --- together with the fact that we know Aristotle's 
analysis only from
lecture notes taken by his students --- explains methinks their 
somewhat muddled-up
phrasing and sometimes cumbersome argumentative development.\\

The Runner (The Dichotomy) [{\small DK} {\scriptsize{29A (25)}}]

\begin{quotation} [Arist., {\em Phys.}, Z 9, 239b(11)] {\em The {\bf 
{\em first}}
[argument] is the one which declares movement to be impossible 
because, however near the
mobile is to any given point, it will always have to cover the half, 
and then the half of
that, end so without end  before it gets to the goal. (...) Hence Zeno's
argument makes a false assumption in asserting that it is impossible 
to pass over an
infinity or to touch one by one infinitely many in a finite time. For 
there are two
senses in which length and time and the continuum in general are said to be
unbounded: with respect to partition and with respect to the extended parts.
Therefore it cannot be assumed possible to touch an infinite quantity 
of things [i.e.,
parts] in a finite time, though this can be assumed for partition, 
because time itself is
infinite in the same way.}
\end{quotation}

\noindent $\rhd$ {\sc Motion cannot take a start}.--- Aristotle reads Zeno's
argument as a {\em reductio}: whether you have to start or arrive, 
seen from whatever
distance you are at, the goal to start from or to arrive at will be 
unattainable.
Interestingly  enough, in the `to start' variant this version of the 
paradox rises
problems even to contemporary limit-based approaches. Apparently the
transition of nothing into something --- from standstill to motion --- is more
problematic than its opposite, from something into nothing, which in 
itself should arouse
suspicion.\footnote{The embarrassment is plain in e.g. Clark's 
commentary on the
``regressive'' formulation of the paradox, see M. Clark, {\em o.c.}, 
p. 160.} The reason
obvioulsy is that the formalisation of the solution of the received 
view --- finite
limits to infinite Cauchy sequences do exist --- is not applicable in 
the `to start'
case. This becomes understandable when we realise that the runner has 
to traverse
$2^{\nn}$ distances to make even his first step! Aristotle correctly interprets
Zeno's intentions, insofar as he admits that two kinds of infinity 
are involved: one with
respect to the --- countable --- number of partitions, one with 
repect to the ---
uncountable --- number of parts. It is well worth the effort to lay 
bare Aristotle's
approach in some detail. First he imputes Zeno with a false 
assumption, viz., ``that it
is impossible to pass over an infinity or to touch one by one 
infinitely many in a finite
time'', while we know that 1) Zeno nowhere mentions time, and 2) this is not an
assumption, but a distorted version of the conclusion reached to the 
argument on
plurality gone before. Thus Aristotle's strategy involves two steps: 
to introduce `time'
as an underlying hypothesis and to discard the plurality arguments, 
although he certainly
was aware of them. Even more, this same pattern turns up against all 
arguments on motion.
Why? The first and major reason obviously is that Aristotle does not 
want to expose, but
to kill off the paradoxes. This stance is exemplified in the basic 
axiom shoring up both
his metaphysics and his logic, the {\em Principium Contradictionis} 
or contradiction
principle (PC): {\em it is not admissible}\/ that something is and is 
not in any sense at
the same place at the same time [{\em Met,}, $\gamma$ 3, 1005b(19-26); B 2,
996b(30)].\footnote{LOEB$_2$. For an in depth discussion, sse J. 
Lukasiewicz, {\em
\"{U}ber den Satz des Widerspruchs bei Aristoteles}, trans. J.M. 
Bochenski, Georg Omls
Verlag, Hildesheim etc., 1993.} For Aristotle paradoxes are a {\em 
problem} most urgently
in need for a solution. Secondly, that solution has to take a 
specific form. This follows
from his criticism of the solution proposed by Plato, who aimed at 
dismantling Zeno's
plurality-argument (in the form given to it according to him by Zeno's master,
Parmenides).\footnote{K. Verelst and B. Coecke, ``Early Greek Thought 
and Perspectives
for the Interpretation of Quantum Mechanics: Preliminaries for an 
Ontological Approach'',
in: {\em Metadebates on Science. The Blue Book of Einstein meets 
Magritte}, G.C.
Cornelis, S. Smets and J.-P. Van Bendegem, Kluwer Academic Press, 
Dordrecht, 1999, pp.
163-196.} Aristotle's reasoning apparently is that if one were to 
avoid the flaws in
Plato's system, one were to avoid the paradoxes of plurality as well, 
and turn instead to
the paradoxes of motion. This, however, is impossible without 
mutilating them. To attain
his goal, the Stagirite proceeds in a most subtle way. Indeed it is 
impossible to
``touch'' stepwise the infinity of parts because the infinities 
involved are different.
But counting is a method of time-measurement --- think of the metre 
in poetry or music ---
therefore the divisibility of time will be `parallel' to that of 
partition. So instead of
looking at partitions and parts, let us look at partitions {\em over 
time}. Aristotle
most cleverly shifts our attention from (discrete) counting and 
(continuous) extension to
merely counting (steps) and counting (time). No wonder that the 
problem disappears! This
is why he says that Zeno's argument fails because time itself is 
infinite ``in the same
way''.  Moreover, however long we {\em count}, e.g., by 
``touching'', the number of
``touches'' remains finite. You only {\em approach}\/ infinity, you 
never reach it. Thus
one never trespasses the PC by making the inevitable cardinal jump 
that is so proper to
Zeno's procedure. Once you introduced time, you can postpone that 
fatal moment as long as
you wish. This is Aristotle's famous {\bf
\em potential} \/infinity. Even if a stretch of time itself is finite,
it consists of a potential infinity of very very small but 
nevertheless {\em finite}
parts: the faster you count, the smaller the parts. The paradox 
disappears because the
very large but finite number of parts of any extended body at 
whatever finitely remote
moment coincides to this same potential infinity, which is exactly 
what countability
means. You always use commensurable quantities, or in modern terms, 
you relate rationals
to rationals. This corresponds to what we continue to do by using the notion of
a mathematical limit. You can maintain you reached the endpoint (of 
the racecourse, say)
because the gap that separates you from it becomes so small that you 
can neglect it. But
this explains as well why no explicit {\em construction}\/ of the 
irrational numbers as
such --- and not merely some arbitrary objects considered azquivalent 
to them --- is
available. G\"{o}del remarks with sore precision: {\em It is 
demonstrable that the
formalism of classical mathematics does not satisfy the vicious 
circle principle in its
first form, since the axioms imply the existence of real numbers 
definable in this
formalism only by reference to all real numbers.}\footnote{K. 
G\"{o}del, ``Russell's
Mathematical Logic'', {\em BPPM}, p. 455.} Or more caustically: {\em 
The calculus
presupposes the calculus.}\footnote{L. Wittgenstein, {\em 
Philosophical Remarks},
Blackwell, Oxford, 1975, p. 130.} With Zeno's paradoxical
construction one does not suffer from this kind of circularity.\\

The Achilles [{\small DK} {\scriptsize{29A (26)}}]

\begin{quotation} [Arist., {\em Phys.}, Z 9, 239b(14)] {\em The {\bf 
{\em second}}
[argument] is the so-called Achilles. This is that the slowest runner 
will never be
overtaken by the swiftest, since the pursuer must first reach the 
point from which the
pursued started, and so the slower must always be ahead. This 
argument is essentially the
same as that depending on dichotomy, but differs from it in that the added
lengths are not divided into halves.}
\end{quotation}

\noindent $\rhd$ {\sc Motion cannot come to its end}.--- Here the 
case is simple, for
Aristotle comments: {\em This argument [the Achilles] is the
same as the former which depends on dichotomia} [{\em Phys.} Z, 9, 
239b (20-21)].
The case is slightly more complicated by the fact that both the 
moving body and the goal
to attain are themselves in motion, but the complication is not substantial, as
Aristotle poinst out himself: it merely implies that the distances to 
cross will not
decrease symmetrically. He mentions {\em
$\delta\iota\chi o\tau o\mu\epsilon\csup{\iota}{-2}\nu$} [{\em Phys,}, Z 3,
239b(19)], symmetric two-division (in his treatment implicitly 
oriented and stepwise
decreasing), but takes care to make clear that even if another number 
of division is used,
this would not make any real difference, something we already know 
from Plato's {\em
Philebus}. The received view presents us Zeno's argumentation as flawed by an
elementary mathematical error, due to a lack of mathematical sophistication. In
accordance with Aristotle's distorted rendering of Zeno's argument, it is
presented as a potentially infinite sequence decreasing geometrically:
$\Sigma{\frac{1}{2^n}}$, the sum of which can very well have a finite total,
because the underlying sequence converges to its finite 
Cauchy-limit.\footnote{R.
Courant and F. John, {\em o.c.}, vol. I, p. 70 sq.} As we already 
mentioned Vlastos,
given his direct acquaintance with the sources of ancient Greek thought at
unease with this modern self-sufficiency, looks for other
explanations,\footnote{VLAS, p. 234.} and proposes alternative 
interpretations based on
the notion of `supertask', developed in the fifthies by Thomson and
Black.\footnote{M. Black, `Achilles and the Tortoise', {\em 
Analysis}, {\bf XI},
1950, pp. 91-101; J. Thomson, `Tasks
and Super-Tasks', {\em Analysis}, {\bf XV}, 1954, pp. 1-13.} A 
supertask requires the
execution of a countable infinity of acts in a finite stretch of 
time. Thomson and
Black argue that such actions are --- under specific circumstances 
--- carried out in
reality, and that they can be used to explain Zeno's motion 
paradoxes. But although we
ever only make finitely many steps, even if we could make countably 
many, the stretches
to croos would be uncountable in number. Thus a Zenonian supertask 
properly speaking would
require an infinity of acts in no time!\\

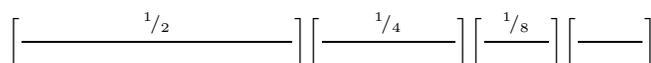
\begin{figure}
\[
\vcenter{\xymatrix@1{{\Bigl[} \ar@{-}[rrrr]^-{^{1}\!/_{2}} &&&& {\Bigr]\Bigl[}
\ar@{-}[rr]^-{^{1}\!/_{4}} && {\Bigr]\Bigl[} 
\ar@{-}[r]^-{^{1}\!/_{8}} &  {\Bigr]\Bigl[}
\ar@{-}[r] & {\Bigr]}}}
\]
\caption{The Received view on the motion paradoxes}
\end{figure}

The Stadium [{\small DK} {\scriptsize{29A (28)}}]

\begin{quotation} [Arist., {\em Phys.}, Z 9, 239b(33)] {\em (a) The 
{\bf {\em fourth}} is
the one about the two rows of equal bodies which move past each other 
in a stadium with
equal velocities in opposite directions, the one row originally 
stretching from the goal [to the
middle-point], the other from the middle-point [to the starting 
point]. This, he thinks,
involves th conclusion that half a given time is equal to its double. 
The fallacy lies in
assuming that a body takes an equal time to pass with equal velocity 
a body that is in
motion and a body of equal size at rest (...) .}
\end{quotation}

\noindent $\rhd$ {\sc The first two arguments combined for extended
bodies}.--- This argument, traditionally known as the fourth argument 
logically is the
third, because it simply combines the former two. What happens when 
two measuring rods
--- our model for Zeno's extended bodies, and one used in this case 
by Aristotle as
well --- pass each other at constant velocity in opposite directions? 
So we now not only
consider the relation rest/motion, but motion/motion as well. Of 
course there still is
the fixed measuring rod with respect to which division through and 
through takes place
(and which creates the impression of time and direction): the floor 
of the stadium. {\em
And this for all generations at once.} The Received View here is that 
Zeno did not
understand the (Galilean) relativity governing the motions of bodies 
in inertially moving
frames of reference, as in the case of two cars crossing each other 
with equal speed on a
high way: {\em The unanimous verdict on Zeno is that he was 
hopelessly confused about
relative velocity in this paradox}.\footnote{{\em STF},
http://plato.stanford.edu/archives/sum2004/entries/paradox-zeno/.} 
But in Zeno's
description, every part at every moment faces its doubling by 
division, whether it be in
comparision to a stable measuring rod, or a rod passing by.  The 
problem arises from the
fact that, because of simultaneous through-and-through division, ``to 
double'' here
involves a transition from ordinal to cardinal, from countable to 
uncountable, from
potential to actual infinity. The infamous ``doubling of the times'' 
only takes into
account the potential, stepwise part of the argument. For of course, 
every body, while
being a continuum, ``touches'' (counts) the other one everywhere when it passes
(measures). It remains just the same cardinal problem. Their speed 
proportional to each
other does not change anything to this fact, analoguous to what we 
saw with the Achilles:
they are at every moment passing each other at infinitely many parts, 
which, by facing
each other's unlimited division, count each other's uncountability.

\SelectTips{cm}{10} \UseTips
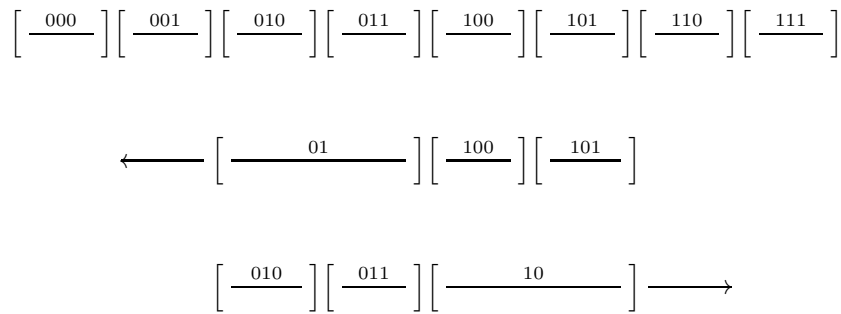
\begin{figure}[h]
\[ {\xymatrix{{\Bigl[} \ar@{-}[r]^-{000} & {\Bigr]\Bigl[}
\ar@{-}[r]^-{001} &  {\Bigr]\Bigl[} \ar@{-}[r]^-{010} & {\Bigr]\Bigl[}
\ar@{-}[r]^-{011} &  {\Bigr]\Bigl[} \ar@{-}[r]^-{100} & {\Bigr]\Bigl[}
\ar@{-}[r]^-{101} &  {\Bigr]\Bigl[} \ar@{-}[r]^-{110} & {\Bigr]\Bigl[}
\ar@{-}[r]^-{111} & {\Bigr]} &\\
&& \ar[l] {\Bigr[}\ar@{-}[rr]^-{01} && {\Bigr]\Bigl[}
\ar@{-}[r]^-{100} & {\Bigr]\Bigl[}
\ar@{-}[r]^-{101} & {\Bigr]} & \\
&& {\Bigl[} \ar@{-}[r]^-{010} & {\Bigr]\Bigl[}
\ar@{-}[r]^-{011} & {\Bigr]\Bigl[}
\ar@{-}[rr]^-{10} && {\Bigr]} \ar[r] &}}
\]
\caption{\label{Figuur-Bewegingsparadoxen}The Stadium}
\end{figure}

\noindent Moreover, the number of partitions -- steps --- involved 
here really is
(countably) infinite. Aristotle's dictum equally applies: {\em For 
there are two senses in
which a distance or a period (or indeed any continuum) may be regarded as
unbounded, viz., with respect to partition and with respect to the 
parts.} Let us stress
once more that Zeno does not imply that motion does not exist, only that it is
paradoxical.  Graphical representations that do not take this fact 
into account do not
account for Zeno's third argument, as is the case with all drawings 
based on Alexander,
as given in Kirk {\em et al.}.\footnote{KRS, p. 274.}\\

The Arrow [{\small DK} {\scriptsize{29A (27)}}]

\begin{quotation} [Arist., {\em Phys.}, Z 9, 239b(30)] {\em The {\bf 
{\em Third}} is that
just given above, that the flying arrow is at rest. This conclusion 
follows from the
assumption that time is composed of instants; for if this is not 
granted the conclusion
cannot be inferred.}
\end{quotation}

\noindent $\rhd$ {\sc Motion is self-contradictory}.---  The Arrow 
radicalises the
reasoning by combining the first two arguments pointwise, so that the 
contradiction
plainly arises. Indeed, even in this case the doubling of parts 
occurs, as again the
Stagirite notes. He does nevertheless not credit Zeno with
this insight. Probably because in this variant the paradox does not 
leave any room for
anything timelike to be smuggled in, Aristotle has to take recourse 
to praemiss of
the parallellism of divisibility of space and time, which is why he 
introduces as an
explicit assumption instants conceived as 
`time-atoms'\footnote{Developed aftewards by
Diodorus Cronus. See R. Sorabji, {\em Time, Creation and the 
Continuum. Theories in
antiquity and the early middle ages}, Duckworth, London, 1983, p. 17 
sq.}, ``for
otherwise the conclusion will not follow'' --- after you have 
discarded actual infinity,
that is. He then rejects these time-atoms, and proposes his potential 
divisibility as a
more apt solution [{\em Phys.}, Z 9, 239a(20-24)]. But of course chronons, like
atoms, are non-Zenonian. This is another nice example of Aristotle's 
general neutralising
strategy with respect to Zeno's paradoxes: to introduce a seemingly 
self-evident
hypothesis on Zeno's behalf, such that his own principle of contradiction can
subsequently be applied succesfully. When parts considered are of the 
{\em megala} type,
one can still be impressed by seeing the motion that takes place. 
When looked at it from
the point of view of the {\em mikra} the paradox becomes unescapable, 
for one cannot see
motion over an unextended ``distance'' in the unextended now. These 
parts in effect
cannot be further divided, which is why the atomistic point of view 
seems to fit in
naturally. But one then forgets an essential thing: the arrow is a 
finite object
consisting of megala {\em and} mikra, which nevertheless flies. The 
last two arguments
show that, whether we consider the {\em megala} or the {\em mikra}, 
the paradox remains
the same. This is where the reading of the last two arguments as a 
dilemma stems from.
Indeed Zeno's argumentation becomes here somewhat of a mocking 
variant of the dilemma
imposed on him by later times. Considered as an absurdity, it seldomly is
discussed with the zeal devoted to the other paradoxes of motion. But 
once the true nature of Zeno's paradoxes is
assimilated, this last argument reveals itself as the contrary of how 
it is generally
perceived: a clear and incontestable exposure of the paradoxical 
nature of motion and
change, and not an incomprehensible enigma. If you let motion,
conceived of as covering all systematically smaller extended parts of 
a line by counting
the uncountable in every single part, come to an end by mentally 
letting the extension of
the parts decrease to nought, then division `comes to an end' too, 
and the only thing that
remains is the naked paradox. This explains why this paradox in the 
literature has been
considered as the most enigmatic one, while it actually only sums up 
Zeno's conclusion in
a concise way.\\

{\sc Geometrical implications}.\footnote{On this subject I was helped 
enormously in
developing the ideas sketched here by some clarifying discussions 
with P. Cara (VUB)
and F. Buekenhout (ULB).}-- To conclude, I will give a short sketch 
of some geometrical
consequences --- to be worked out in more detail in a subsequent 
paper --- which throw
light on Zeno's paradoxes from a more physical perspective, and which 
allow to bring the
Received View on the paradoxes of motion into the picture again. At 
first an observation
that serves as a guideline. We will work in the spirit of the 
notorious {\em Erlanger
Programm}, formulated by F. Klein in 1872.\footnote{F. Klein, ``Vergleichende
Betrachtungen \"{u}ber neuere geometrische Forschungen'', {\em 
Mathematische Annalen},
{\bf 43}, p. 63 sq., 1893.} Instead of focusing on geometric objects 
{\em per se}, one
studies objects that stay invariant under the action of a group of 
transformations. Now
the algebraic expression for Nieuwentijt infinitesimals is a first 
order polynome, the
equation of an Euclidean straight line. Such a polynome $a + \eta b$ 
can also be regarded
as a point $(a, b)$ \/in a two-dimensional space. The usefulness of 
this representation is
clear from the example of the complex numbers $g = a + ib$\,, with 
\,$i^2 = -1$ for \,
$\eta = i \,:$; their geometrical representation is the complex face 
of the Euclidean
plane with its accompanying static, Euclidean geometry. Let us call 
$a + \eta b$ a two
dimensional number\footnote{A suggestion by K. Lefever.} and ask 
whether for other values
of $\eta$ a geometrical representation exists.\footnote{Emch} Is 
there a likely candidate
to fulfill such a role for our Nieuwentijt infinitesimals $g = a + 
hb$\,; $h^2 = 0$, with
\, $\eta = h
\;$? Although less well known, there indeed is: {\bf Galilean 
geometry}.\footnote{I.M.
Yaglom, {\em A Simple Non-Euclidean Geometry and its Physical Basis}, 
Springer-Verlag, N.
Y., 1979.} The geometry in which the motions of bodies in classical 
mechanics take place
is non-Euclidean! We see that the first order polynome describing 
each Nieuwentijt number
is the equation of a uniform linear motion in Euclidean space. The 
transition from static
to kinematic is not so innocent as it would seem, for it is by 
introducing {\em time} that
one captures motion; in which case limit-like approximations become 
possible, and the
traditional view on PM can be recovered. Euclidean geometry gives us 
merely the situation
of a particle instantenously. Now if we take the hyperreal part
$bh$ as representing the velocity-component of the Galilean 
transformation $x' = x + vt$
\,(with $h$\; as Galilean time $t$\,), we see that only on the {\em 
loci} of the {\em
dubbels} the hyperreal monads have the power to generate a space-time 
worldline: the idea
is that there some {\em action} is involved. We saw also that we 
nevertheless will get
back the Euclidean plane. When the
$bh$-coordinate on the time axis of a hypereal monad will differ
from zero, transition from Euclidean to Galilean geometry takes 
place, i.e., the thing
sitting on these coordinates is set into motion. But this way of 
building Galilean
geometry seems awkward, since it implies that points can transform 
into straight lines.
However, a geometry does exist in which certain points can explode 
into ``higher order''
points under specific conditions, i.e., lines, or even curves. This 
is {\bf Cremona
geometry}\footnote{F. Enriques, {\em Lezioni sulla teoria geometrica 
delle equazioni e
delle funzioni algebriche}, N. Zanichelli Editore, Bologna, 1915.}, 
the geometry of
birational transformations of the plane. Let us summarise its basic 
tenets and show by
analogy of argumentation that it is reasonable to expect it to be the 
appropriate
geometrical description of the transformations implicit in Zeno's approach. \\

Let a {\em birational transformation} \/$\varphi$ be a birational 
function of the
coordinates of $x$, whereby \,$\varphi : V \pijlr V' : x \pijlr x'$. 
The function
$\varphi$ is not defined everywhere for at certain points the 
denominator will become
zero and singularities will arise. Thence points do not always have 
an image under
transformation, so they cannot be invariant. Those who are not are 
called singular
points; we can think of them as `black' or `invisible'. We are here 
in the realm of
algebraic geometry. In order to get rid of our divisions by zero, we 
work preferably in
the completed, projective plane. In that case the transformation can 
be written as
homogeneous polynomials of the coordinates of
$x$.\footnote{In which case they can be simply rational. J. Tits, 
{\em The Cremona Plane},
Lecture notes by H. Van Maldeghem and F. Buekenhout, VUB-ULB, 1999.} 
Singular points are
characterised by higher order `points', or approximations in their 
infinitesimal
neighbourhoods. Such points can then be taken as transformed into 
curves by division
through zero. They are represented by polynomes of order $n$, in the 
same way analysis
approximates functions by means of Taylor expansions.\footnote{F. 
Enriques, {\em o.c.},
vol. II, p. 327 sq.} In our case all approximations will be first 
order polynomes
representing Nieuwentijt infinitesimals at every {\em locus} where a double
$[d_i, d'_i]$ exists, thus transforming singular points into straight 
lines. Physically
speaking this comes down to derivation, another way to ascertain that 
the traditional,
time-dependent perspective on the paradoxes of motion is included in 
our approach. Owing
to Sierpi\'{n}ski's theorem mentioned above, our completed plane 
would be a special
instance of the projective plane over the division ring of nilpotent numbers
$\bb$; we will label it $P_2(\bb)$ (this remains to be rigorously 
shown).\footnote{F.
Buekenhout and P. Cameron, ``Projective and Affine Geometry over 
Division Rings'', in F.
Buekenhout (ed.), {\em Handbook of Incidence Geometry}, Elsevier, Amsterdam,
1995, pp. 27-62.} Cremona Geometry was axiomatised by J. Tits in the 
context of incidence
geometry, in which geometrical properties are expressed as 
symmetrical relations of
intersection and inclusion.\footnote{J. Tits: {\em Th\'{e}orie des groupes},
Coll\`{e}ge de France, R\'{e}sum\'{e} des cours et travaux, 
1998-1999.} In every ``black
point'' Tits defines a  tree of approximations which resembles
$({\cal Z}, \order)$. Far from being empty, the singular point 
appears to be a highly
structured entity! These trees are the foundation on which the 
Cremona plane can be
constructed as a {\em building}, incidence geometrically speaking. 
They constitute thin
subgeometries, and are called {\em apartments}.\footnote{In Tits's 
original vocabulary
these objects were called ``squelettes'' and ``ossuaires''!} What is 
the precise nature
of these apartments? Given the strong analogy between the structures in Tits's
axiomatisation and the semi-lattices arising from our Zeno-approach, 
F. Buekenhout
proposed the following theorem:

\begin{quotation}\noindent {\em The `thin' Zeno-plane (the small 
Zeno-semilattice ${\cal
Z}$) gives the [aequivalence class of] apartments in the building 
constituted by the
Cremona plane.}\end{quotation}

\noindent It of course remains to be demonstrated that this theorem 
does indeed establish
the desired link. When this works it means that the ``Zeno-line'' is 
identical to the
Cremona line. In that case one could use a kind of ``Zeno 
microscope'' to elucidate the
internal structure of the Nieuwentijt infinitesimal, which, far from being
structureless, repeats fractal-wise the Zenonian tree in its own fine 
structure. This
could be the first step towards an understanding of why the `cardinal 
jump' so crucial to
Zeno's paradoxical procedure comes about. The question then remains 
to be answered what
the birational transformations are under which Nieuwentijt 
infinitesimals are invariant.
Settling this question would bring us back to Galilean geometry, 
i.e., to classical
mechanics, which really is where Zeno's paradoxes of motion belong.

\bigskip

{\sc Acknowledgments}.-- I want to express my gratitude to the people 
of the Department
of Mathematics at the Vrije Universiteit Brussel, for not only 
tolerating but even
actively encouraging a (l)on(e)ly philosopher working among them and 
foraying their
noble discipline for what at first instance could only have appeared as wild
speculations. I especially thank Tim Vanderlinden and Philippe Cara, 
who not only
critically followed my attempts at precise mathematical formulation, 
but even contributed
substantially to the development of my ideas. Furthermore I thank 
Rudger Kieboom,
Tomas Everaert and Mark Sioen (dept. Math., VUB), Wendy Lowen 
(Instit. de Math. de
Jussieu, Paris VII), as well as Roland Hinnion and Francis Buekenhout 
(dept. Math.,
Universit\'{e} Libre de Bruxelles). Bob Coecke (Comlab, Oxford University) is a
longstanding support, and introduced me to the secrets of Domain 
Theory. To Rudolf De
smet and Wilfried Van Rengen (Classical Studies, VUB) I am indebted 
for competent advise
with respect to philological issues. A final word of thanks to Koen 
Lefever, whose
intellectual influence is less specific, but acts as a light-house to 
my ongoing
philosophical quest.

\bigskip

{\sc A list of bibliographical sigla}.-- Abbreviations used 
throughout the text are
assembled here. For further bibliographical information the reader is 
referred to the
footnotes.\\
\par\noindent
{\em DK} :  H. Diels and W. Kranz, {\em Fragmente der Vorsokratiker}, 
erster Band,
Weidmann, Dublin, Z\"{u}rich, 1951 [1996].
\par\noindent
{\em BPPM} : P. Benacerraf \& H. Putnam, {\em Philosophy of 
Mathematics}, Cambridge
University Press, Cambridge, 1964 [1991].
\par\noindent
{\em FTG} : J. Van Heijenoort, {\em From Frege to G\"{o}del}, Harvard 
University Press,
Cambridge Mass., 1967.
\par\noindent
{\em GA} : G. Cantor, {\em Gesammelte Abhandlungen}, E. Zermelo Ed., 
Georg Olms Verlag,
Hildesheim, 1932 [1962].
\par\noindent
{\em LEE} :  H.D.P. Lee, {\em Zeno of Elea. A Text, with
Translation and Notes}, Cambridge University Press, Cambridge, 1936.
\par\noindent
{\em KRS} : G.S. Kirk, J.E. Raven and M. Schofield, {\em The 
Presocratic Philosophers. A
critical History with a Selection of Texts}, Cambridge Univ. Press, Cambridge,
1957 [1983].
\par\noindent
{\em LOEB$_1$} : Aristotle, {\em Metaphysics}, Books I-IX, transl. H. 
Tredennick, Harvard
University Press, Harvard, 1933 [1996].
\par\noindent
{\em LOEB$_2$} : Aristotle, {\em Physics}, Books I-VIII, transl. P.H. 
Wicksteed and
F.M. Cornford, Harvard University Press, Cambridge, Mass., 1933 [1995].
\par\noindent
{\em LOEB$_3$} : Plato, Vol. I, trans. H.N. Fowler, Harvard 
University Press, Cambridge,
Mass., 1914 [1999].
\par\noindent
{\em LOEB$_4$} : Plato, Vol. VIII, trans. H.N. Fowler and W.R.M. 
Lamb, Harvard University
Press, Cambridge, Mass., 1925 [1999].
\par\noindent
{\em SHYP} :  H. Poincar\'{e}, {\em La science et l'hypoth\`{e}se}, 
Flammarion, Paris,
1902 [1968].
\par\noindent
{\em STF} : E.N. Zalta (ed.), The Stanford Encyclopedia of Philosophy 
(on the Internet).
\par\noindent
{\em SVAL}, H. Poincar\'{e}, {\em La valeur de la science}, 
Flammarion, Paris, 1970
[1905].
\par\noindent
{\em VLAS} : G. Vlastos, {\em Studies in Greek Philosophy. Vol. I: The
Presocratics}, Princeton University Press, Princeton, 1993.

\end{document}